\documentclass[12pt]{article}

\UseRawInputEncoding

\usepackage{mathtools,amsthm,amssymb,color,cases}

\usepackage{stackengine}

\usepackage{mathrsfs}
\usepackage[abbrev]{amsrefs}
\usepackage{tikz}
\usepackage{enumitem}
\usepackage{setspace}

\usepackage{latexsym}
\usepackage{mathrsfs}
\usepackage{comment}
\usepackage{microtype}

\usepackage{indentfirst}

\usepackage{geometry}
\usepackage{etoolbox}
\AtEndEnvironment{proof}{\setcounter{proofpart}{0}}

 \newtheoremstyle{mystyle}
    {}
    {}
    {\normalfont}
    {}
    {\bfseries}
    {}
    { }
    {}
  \theoremstyle{mystyle}     

\allowdisplaybreaks[1]
\theoremstyle{theorem}
\newtheorem{thm}{Theorem}[section]

\newtheorem{lem}[thm]{Lemma}
\newtheorem*{rem*}{Remark}

\newtheorem{prop}[thm]{Proposition}

\theoremstyle{definition}
\newtheorem{definition}{Definition}

\newtheorem{rem}[thm]{Remark}

\newtheoremstyle{part}{}{}{\normalfont}{}{\itshape}{.}{.5em}{}
\theoremstyle{part}

\newcommand\blfootnote[1]{%
  \begingroup
  \renewcommand\thefootnote{}\footnote{#1}%
  \addtocounter{footnote}{-1}%
  \endgroup
}

\numberwithin{equation}{section}
\numberwithin{thm}{section}

\newenvironment{equ*}{
    \begin{equation*}
}{
    \end{equation*}
}

\newenvironment{pf}
   {{\noindent \bf Proof.}}{\hfill \qed}

\newtheoremstyle{part}{}{}{\normalfont}{}{\itshape}{.}{.5em}{}
\theoremstyle{part}

\newcommand{\ddj}{\dot{\Delta}_j}

\newcommand{\pt}{\partial}

\newcommand{\F}{\mathcal {F}}  
\newcommand{\G}{\mathcal {G}}  
\renewcommand{\S}{\mathcal {S}}
\newcommand{\Z}{\mathbb{Z}}

\newcommand{\R}{\mathbb R}
\newcommand{\Nat}{\mathbb N}

\newcommand{\intr}{\int_{\R^3}}

\newcommand{\al}{\alpha}

\newcommand{\x}{\xi}

\renewcommand{\t}{\theta}

\newcommand{\N}{\nabla }

\renewcommand{\div}{{\rm {div}}}

\newcommand{\supp}{\text{ supp }}

\renewcommand{\P}{\mathcal{P}}
\newcommand{\Q}{\mathcal{Q}}

\newcommand*{\dd}{\mathop{}\!\mathrm{d}}

\newcommand{\ha}{\hat{a}}
\newcommand{\hv}{\hat{v}}
\newcommand{\hu}{\hat{u}}
\newcommand{\hp}{\hat{\phi}}
\newcommand{\hps}{\hat{\Psi}}

\newcommand{\db}{\dot{B}}
\newcommand{\fdp}{ \frac{d}{p} }

\DeclarePairedDelimiter{\norm}{\lVert}{\rVert}
\DeclarePairedDelimiter{\bignorm}{\bigg\|}{\bigg\|}

\DeclarePairedDelimiter{\floor}{\lfloor}{\rfloor}

\DeclarePairedDelimiter{\bigc}{\big{(}}{\big{)}}
\DeclarePairedDelimiter{\Bigc}{\Big{(}}{\Big{)}}
\DeclarePairedDelimiter{\Biggc}{\Bigg{(}}{\Bigg{)}}
\DeclarePairedDelimiter{\biggc}{\bigg{(}}{\bigg{)}}

\DeclarePairedDelimiter{\Bigf}{\Big{[}}{\Big{]}}



\renewcommand{\db}{\dot{B}}
\renewcommand{\fdp}{ \frac{d}{p} }

\newcommand{\NL}{
    \int^t_0
    e^{(t-s) M(D) } 
    \begin{bmatrix}
    f(s) \\
    h(s)
    \end{bmatrix}
    \dd s
}

\renewcommand{\G}{
\mathcal{G}
}

\begin{document}

\begin{center}
    {\bf \large Optimal decay estimates for the radially symmetric compressible Navier-Stokes equations} \\
    
    \phantom{}
    
    Tsukasa Iwabuchi* \quad D\'aith\'i \'O hAodha** \\
    
     \phantom{}
    
    Mathematical Institute, Tohoku University,
    
    980-0845
\end{center}

\blfootnote{Email: *t-iwabuchi@tohoku.ac.jp, **david.declan.hughes.p6@dc.tohoku.ac.jp}

\begin{center}
\begin{minipage}{135mm}
\footnotesize
{\sc Abstract. }
We examine the large-time behaviour of solutions to the compressible Navier-Stokes equations under the assumption of radial symmetry. In particular, we calculate a fast time-decay estimate of the norm of the nonlinear part of the solution.
This allows us to obtain a bound from below for the time-decay of the solution in $L^\infty$, proving that our decay estimate in that space is sharp.
The decay rate is the same as that of the linear problem for curl-free flow.
We also obtain an estimate for a scalar system related to curl-free solutions to the compressible Navier-Stokes equations in a weighted Lebesgue space.
\end{minipage}
\end{center}

\section{Introduction}

In this paper, we consider the barotropic compressible Navier-Stokes system
\begin{align} \label{CNSo}
    \begin{cases}
    \pt_t \rho + \div (\rho u) = 0, & \text{ in } (0,\infty) \times\R^3, \\
    \pt_t (\rho u) + \div (\rho u \otimes u) - \div(2 \mu D(u) + \lambda \, \div (u) \text{Id}) + \nabla p = 0, & \text{ in } (0,\infty) \times\R^3, \\
    (\rho, u) |_{t=0} = (\rho_0, u_0), & \text{ in } \R^3,
    \end{cases}
\end{align}
where $\rho : [0,\infty) \times \R^3 \to [0,\infty),$ and $u: [0,\infty) \times \R^3 \to \R^3$ are unknown functions, representing the density and velocity of a fluid, respectively. 
$p: [0,\infty) \times \R^3 \to \R$ is the pressure in the fluid, and the barotropic assumption gives us $p \coloneqq P(\rho)$, for some smooth function $P(\cdot)$.
$\mu,\lambda$ are viscosity coefficients, taken such that
\begin{align*}
    \mu > 0, \quad 2\mu + \lambda > 0.
\end{align*}
We define the deformation tensor
\begin{align*}
    D(u) \coloneqq \frac{1}{2} \Big{(} Du + Du^T \Big{)}.
\end{align*}


In this paper, we will obtain time-decay estimates of solutions to the radially symmetric case of the above problem. 
Before we introduce our main result, we discuss a few  well-known results concerning \eqref{CNSo}.
Matsumura-Nishida showed in~\cites{matsumura-nishida1979, matsumura-nishida1980} that \eqref{CNSo} has global solutions when equipped with data $(\rho_0, u_0)$ that is a small perturbation in $L^1 \cap H^3$ of $(\bar{\rho}, 0)$ for any positive constant $\bar{\rho}$, and proved the following decay result
\begin{align*}
    \bignorm{ 
    \begin{bmatrix}
        \rho(t) - \bar{\rho} \\
        u(t)
    \end{bmatrix}  
    }_{2}
    \leq C ({t+1})^{-3/4}.
\end{align*}
This is the decay rate of the solution to the heat equation with initial data in $L^1$.
Ponce then extended these results to other $L^p$ norms. In particular, for $p \in [2,\infty]$, $k\in\{0,1,2\}$, and dimension $d=2,3$,
\begin{align*}
    \bignorm{ 
    \nabla^k \begin{bmatrix}
        \rho(t) - \bar{\rho} \\
        u(t)
    \end{bmatrix}  
    }_{p}
    \leq C ({t+1})^{-\frac{d}{2}(1-\frac{1}{p}) - \frac{k}{2} }.
\end{align*}

Our results make use of developments in the theory of the compressible Navier-Stokes equations in Besov spaces. Global existence of strong solutions to \eqref{CNSo} for small initial data $(\rho_0, u_0)$ in critical Besov spaces $ \dot{B}^{d/2}_{2,1} \times \dot{B}^{d/2 - 1}_{2,1} $ was first proven by Danchin in~\cite{Danchin2000} and large-time estimates in Besov norms for $p$ close to $2$ were proven by Danchin-Xu in~\cite{Danchin-Xu2017}.
The authors of the present paper proved optimality of decay estimates for the linear problem in Besov spaces in \cite{ohaodha-iwabuchi2023(LinCNS)},
while global existence of solutions in critical weighted Besov spaces were recently proven by the second author in~\cite{ohaodha2023(WBCNS)}.

Our goal in this paper is to obtain an optimal bound of the solution to the system \eqref{CNSo} under the assumption that the initial data $(\rho_0, u_0)$ is radially symmetric.
That is, for all $x\in\R^3$,
\[
\rho_0(x) = \rho_0(|x|), \quad u_0(x) = U_0 ( |x|) \frac{x}{|x|}, 
\]
where $U_0 : [0,\infty) \to \R$.
In particular, we prove a bound from above, in terms of time $t$, of the norm of solutions over space $x$.
By expressing the solution $(\rho, u)$ as the solution to the integral equation (i.e. by considering mild solutions), we shall obtain separate bounds for the linear and nonlinear parts of the solution.
Thanks to the radial symmetry, we show that the nonlinear term decays faster than the linear term.
Then, using the bound from below for the linear term proven in~\cite{ohaodha-iwabuchi2023(LinCNS)}, we can show that the decay rate obtained for the whole solution is sharp.

We make extensive use of several existence and decay results in order to obtain our main theorem. 
First, Hoff-Zumbrun prove the following existence and decay result in \cite{hoff-zumbrun}.
\begin{prop} {\rm(}\cite{hoff-zumbrun}{\rm)}
\label{Hoff-Zumbrun}
    Let $m\coloneqq \rho u$, $m_0 \coloneqq \rho_0 u_0$. Assume that 
    \[
    E \coloneqq 
    \bignorm{  
    \begin{bmatrix}
        a_0 \\
        m_0
    \end{bmatrix}
    }_1
    +
    \bignorm{
    \begin{bmatrix}
        a_0 \\
        m_0
    \end{bmatrix}}_{H^{ 1 + l }}
    \]
    is sufficiently small, where $l\geq3$ is an integer.
Then the Navier-Stokes system \eqref{CNSo} with initial data $\rho_0, u_0$ has a global solution satisfying the following decay estimate 
for any multi-index $\al$ with $|\al| \leq (l-3)/2$:
\begin{align} \label{HZ not curl-free}
        &
        \bignorm{
        D^\al_x \Bigc{
        \begin{bmatrix}
         a(t) \\
         m(t)
        \end{bmatrix}}
        }_p
        \leq
        C(l) E 
    \begin{cases}
    (t+1)^{ - \frac{3}{2} (1-\frac{1}{p})  },
    & 2 \leq p \leq \infty,
    \\
    (t+1)^{ - \frac{3}{2} (1-\frac{1}{p})  - \frac{1}{2} ( 1 - \frac{2}{p} )   - \frac{|\al|}{2}   }, & 1 \leq p < 2,
    \end{cases}
        \\
        &
        \bignorm{
        D^\al_x \Bigc{
        \begin{bmatrix}
         a(t) \\
         m(t) - e^{t\Delta} \P m_0
        \end{bmatrix}}
        }_p
        \leq
        C(l) E (t+1)^{ - \frac{3}{2} (1-\frac{1}{p})  - \frac{1}{2} ( 1 - \frac{2}{p} )   - \frac{|\al|}{2}   },
        \ 2\leq p \leq \infty
        . \label{HZ curl-free}
\end{align}
\end{prop}
Note that, in the norm in inequality \eqref{HZ curl-free}, we are removing the divergence-free part of the linear term of $m$. Thus what remain are the nonlinear term and the curl-free part of the linear term.

Kobayashi-Shibata in \cite{kobayashi-shibata} obtained an estimate for a linearised version of \eqref{CNSo} which separates the solution into high and low frequencies (see Definition \ref{besov spaces} below).
The decay rate in \eqref{HZ not curl-free} is associated with the low-frequency part of solutions, while the high-frequency part decays exponentially with $t$.

In this paper, we will assume that the density approaches $1$ at infinity; and so we are concerned with strong solutions which are small perturbations from a constant state $(\rho,u) \equiv (1,0)$. 
We shall also assume that $\mu, \lambda$ are constant,
and set $a \coloneqq \rho - 1$. 
Our system \eqref{CNSo} can thus be rewritten into the following linearised problem:
\begin{align}
\label{CNS}
    \begin{cases}
    \pt_t a + \div(u) = f & \text{ in } (0,\infty) \times \R^3, \\
    \pt_t u - \mathcal{A} u  + P'(1)\nabla a = g & \text{ in } (0,\infty) \times \R^3, \\
    (a,u)\Big{|}_{t=0} = (a_0,u_0) & \text{ in }\R^3,
    \end{cases}
\end{align}
where the differential operator $\mathcal{A}$ 
is defined by:
\[ \mathcal{A}u \coloneqq \mu \Delta u + (\lambda + \mu) \nabla \div(u), \]
and where the nonlinear terms $f,g$ are defined as follows:
\begin{align*}
    f & \coloneqq -\div(au),
    \\
    g & \coloneqq 
    -u\cdot\N u - \frac{a}{1+a} \mathcal{A} u - \beta(a) \N a,
\end{align*}
with 
\[\beta(a) \coloneqq  \frac{P'(1+a)}{1+a} - P'(1).\]

We make regular use of two results for the problem \eqref{CNS}, both of which use the Besov framework, which we introduce now.
\begin{definition} \label{besov spaces}
Let $\{ \hat{\phi}_j \}_{ j \in\Z}$ be a set of non-negative measurable functions such that 
\begin{enumerate}
    \item $\displaystyle \sum_{ j \in\Z} \hat{\phi}_j (\x) = 1, \text{ for all } \x \in \R^3 \backslash \{0\}$,
    \item $\hat{\phi}_j (\x) = \hat{\phi}_0(2^{-j}\x)$,
    \item $\supp \hat{\phi}_j (\x) \subseteq \{ \x \in \R^3 \ | \ 2^{j-1} \leq |\x| \leq 2^{j+1} \}$.
\end{enumerate}
For a tempered distribution $f$, we write
\[
\dot{\Delta}_j f \coloneqq \F^{-1} [\hat{\phi}_j \hat{f}].
\]
This gives us the \textit{Littlewood-Paley decomposition} of $f$:
\begin{align*}
    f = \sum_{j\in\Z} \ddj f.
\end{align*}
This equality only holds modulo functions whose Fourier transforms are supported at $0$, i.e. polynomials.
To ensure equality in the sense of distributions, we
next let $\dot{S}_j$ denote the low-frequency cutoff function.
That is, for $j\in\Z$,
\begin{align*}
    \dot{S}_j f \coloneqq \bigc{ \chi_{j}(D) + \dot{\Delta}_{j}} f,
\end{align*}
where
\begin{align*}
    \chi_j (D) f \coloneqq \F^{-1} [ \chi(2^{-j} \x) \hat{f} ],
\end{align*}
and $\chi$ is the identity function on $\{ x\in\R^3 \ | \ |x| \leq 1 \}$.
Then we consider the subset $\mathcal{S}'_h$ of tempered distributions $f$ such that
\begin{align*}
    \lim_{j\to -\infty} \norm{ \dot{S}_j f }_{L^\infty} = 0.
\end{align*}
The Besov norm is then defined as follows: for $1 \leq p,q \leq \infty$, and $s \in \R$, we define 
\[
\norm{f}_{\dot{B}^{s}_{p,q}} \coloneqq \Bigc{ \sum_{j \in \Z} 2^{sqj} \norm{\dot{\Delta}_j f}^q_{p} }^{ \frac{1}{q} }.
\]
The set $\dot{B}^{s}_{p,q}$ is defined as the set of functions, $f \in \mathcal{S}'_h$, whose Besov norm is finite. 
Throughout this paper, we will refer to the parameter $s$ as the `\textit{regularity exponent}' and $p$ as the `\textit{Lebesgue exponent}.'

We then define a weighted Besov space as the set 
of functions $f \in \mathcal{S}'_h$ such that the Besov norm of $f$ multiplied by $x_k$ is finite for all $k \in \{1, 2, 3\}$. That is,
\[
\norm{x_k f}_{\dot{B}^{s}_{p,q}} < \infty.
\]
We then call $\norm{x_k f}_{\dot{B}^{s}_{p,q}}$ the weighted Besov norm of $f$.

We also regularly use the following notation for so-called high-frequency and low-frequency norms:
\[
\norm{f}^h_{\dot{B}^{s}_{p,q}} \coloneqq \Bigc{ \sum_{j \geq j_0} 2^{sqj} \norm{\dot{\Delta}_j f}^q_{p} }^{ \frac{1}{q}}, \quad
\norm{f}^l_{\dot{B}^{s}_{p,q}} \coloneqq \Bigc{ \sum_{j \leq j_0} 2^{sqj} \norm{\dot{\Delta}_j f}^q_{p} }^{ \frac{1}{q}},
\]
where $j_0 \in \Z$ is called the frequency cut-off constant.
We also define the high-frequency and low-frequency parts of a function $f$:
\[
f^h \coloneqq \sum_{j \geq j_0} \dot{\Delta}_j f, 
\quad
f^l \coloneqq \sum_{j \leq j_0} \dot{\Delta}_j f.
\]
\end{definition}


The first result for \eqref{CNS} in the Besov framework that we use is due to Danchin-Xu, and gives a global existence and decay result for $(a,u)$ in the critical Besov framework.
For this result, we introduce the function space $X_p$ as the set of all pairs of functions $(a,u)$, where $a : [0,\infty) \times \R^d \to [0,\infty)$ is a scalar function and $u: [0,\infty) \times \R^d \to \R^d$ is a $d$-vector function, satisfying the following:
\begin{align*}
    & (a,u)^l \in \widetilde{C} (\R_{>0}; \db^{\frac{d}{2} - 1}_{2,1}) \cap L^1(\R_{>0}; \db^{\frac{d}{2} + 1}_{2,1}), 
    \\
    & a^h \in \widetilde{C} (\R_{>0}; \db^{\frac{d}{p} }_{p,1}) \cap L^1(\R_{>0}; \db^{\frac{d}{p} }_{p,1}),
    \\
    & u^h \in \widetilde{C}(\R_{>0}; \db^{\frac{d}{p}}_{p,1})
    \cap L^1(\R_{>0}; \db^{\frac{d}{p} + 1}_{p,1})
    .
\end{align*}
$X_p$ is then equipped with the obvious norm corresponding to the strong topologies for the above spaces.
The space $X_p$ is the original `critical space' used for the global existence theorems in \cites{danchin2016, Danchin-Xu2017}.

\begin{prop} {\rm (}\cite{Danchin-Xu2017}{\rm )}\label{dx2017}
    Let $d\geq2$ and $p \in [2, \text{min} \{4, 2d/(d-2)\} ]$, with $p\neq 4$ in the $d=2$ case.
    Assume without loss of generality that $P'(1) = \nu = 1$.
    Then there exists a constant $c = c(p,d, \mu, P)>0$ such that if
    \begin{align*}
        X_{p,0} \coloneqq 
        \norm{ (a_0,u_0) }^l_{ \db^{ \frac{d}{2}- 1}_{2,1} }
        +
        \norm{a_0}^h_{ \db^{\fdp }_{p,1} }
        +
        \norm{u_0}^h_{ \db^{ \fdp  - 1 }_{p,1} }
        \leq 
        c,
    \end{align*}
    then \eqref{CNS} has a unique global-in-time solution $(a,u)$ in $X_p$. 
    Furthermore, there exists a constant $C=C(p,d,\mu, P)>0$ such that
    \begin{align*}
        \norm{(a,u)}_{X_p} \leq C X_{p,0}.
    \end{align*}
    Also, there exists a constant $c_1$ such that if, in addition,
    \begin{align*}
        \norm{ (a_0,u_0) }^l_{ \db^{ -s_0 }_{2,\infty} } 
        \leq 
        c_1,
        \text{ where }
        s_0 \coloneqq d\biggc{ \frac{2}{p} - \frac{1}{2} },
    \end{align*}
    then we have a constant $C_1$ such that for all $t\geq 0$,
    \begin{align*} 
        D_{p,\epsilon}(t) \leq C_1 \biggc{
        \norm{ (a_0,u_0) }^l_{ \db^{ -s_0 }_{2,\infty} } 
        +
        \norm{ (\N a_0,u_0) }^h_{ \db^{\fdp  - 1 }_{p,1} } 
        },
    \end{align*}
    where the norm $D(t)$ is defined by 
    \begin{align*}
        D(t)
        &
        \coloneqq
        \sup_{ s\in [\epsilon -s_0, d/2 + 1] }
        \norm{ 
        \langle 
        \tau
        \rangle^{
        (s_0+s)/2
        } 
        (a,u)
        }^l_{ L^\infty_t \db^s_{2,1} }
        \\
        &
        +
        \norm{
        \langle 
        \tau
        \rangle^{ \fdp  + 1/2 - \epsilon}
        (\N a, u)
        }^h_{ \widetilde{L}^\infty_t \db^{ \fdp  - 1 }_{p,1}  }
        +
        \norm{ \tau \N u }^h_{ \widetilde{L}^\infty_t \db^{ \fdp  }_{p,1} },
    \end{align*}
    with $\epsilon >0$ taken sufficiently small.
\end{prop}

The next result we use is due to the second author and concerns the global existence of solutions in weighted Besov spaces. 
We introduce the solution space $S$ to which our solution will belong as the set of all pairs of functions $(a,u)$, where $a : [0,\infty) \times \R^d \to [0,\infty)$ is a scalar function and $u: [0,\infty) \times \R^d \to \R^d$ is a $d$-vector function, satisfying the following:
\begin{align*}
    & (a,u)^l \in \widetilde{C} (\R_{>0}; \db^{\frac{d}{2} - 1}_{2,1}) \cap L^1(\R_{>0}; \db^{\frac{d}{2} + 1}_{2,1}), 
    \quad
    a^h \in \widetilde{C} (\R_{>0}; \db^{\frac{d}{2} + 1}_{2,1}) \cap L^1(\R_{>0}; \db^{\frac{d}{2} + 1}_{2,1}),
    \\
    & u^h \in \widetilde{C}(\R_{>0}; \db^{\frac{d}{2}}_{2,1})
    \cap L^1(\R_{>0}; \db^{\frac{d}{2} + 2}_{2,1}),
    \\
    & (x_k a, x_k u)^l \in \widetilde{C} (\R_{>0}; \db^{\frac{d}{2}}_{2,1}) 
    \cap L^1_t ( \R_{>0} ; \db^{\frac{d}{2} + 2}_{2,1} ), 
    \quad 
    (x_k a)^h \in \widetilde{C} (\R_{>0}; \db^{\frac{d}{2} +1 }_{2,1})
    \cap L^1_t ( \R_{>0} ; \db^{\frac{d}{2} + 1}_{2,1} ),
    \\
    & 
    (x_k u)^h \in \widetilde{C} (\R_{>0}; \db^{\frac{d}{2}}_{2,1})
    \cap L^1_t ( \R_{>0} ; \db^{\frac{d}{2} + 2}_{2,1} ).
\end{align*}
$S$ is then equipped with the obvious norm corresponding to the strong topologies for the above spaces.
Here, $\widetilde{C} (\R_{>0}; \db^{s}_{p,1}) \coloneqq 
{C} (\R_{>0}; \db^{s}_{p,1})
\cap
\widetilde{L}^\infty(\R_{>0}; \db^{s}_{p,1}),
$ for $s\in\R, \ p \in[1,\infty].$ 
The norm of $\widetilde{L}^\infty( 0,T ; \db^{s}_{p,1})$ for $T>0$ is defined by taking the $L^\infty$-norm over the time interval \textit{before} summing over $j$ for the Besov norm. That is, for all $f \in \widetilde{L}^\infty( 0,T ; \db^{s}_{p,1})$,
\[
\norm{f}_{ \widetilde{L}^\infty( 0,T ; \db^{s}_{p,1}) }
\coloneqq
\sum_{j\in\Z} 2^{sj} \sup_{ t\in(0,T) } \norm{ \ddj f(t) }_{L^p} .
\]
We abbreviate the notation for norms by writing 
\[
\norm{f}_{ \widetilde{L}^\infty_T \db^{s}_{p,1} }
\coloneqq
\norm{f}_{ \widetilde{L}^\infty( 0,T ; \db^{s}_{p,1}) },
\]
and similarly abbreviate other norms over time and space.

\begin{prop} {\rm(}\cite{ohaodha2023(WBCNS)}{\rm)} 
\label{3rd paper main thm}
    Let $d\geq 3$. 
    Assume $P'(1) > 0.$
    Then there exists 
    a frequency cut-off constant $j_0 \in \Z$ and
    a small constant $c = c(d,\mu,P) \in \R$ such that, if $(a_0,u_0)$ 
    satisfy 
    \begin{align*}
        S_{0} \coloneqq &  
        \norm{(a_0,u_0)}^l_{\dot{B}^{\frac{d}{2} - 1}_{2,1}} 
        + \norm{a_0}^h_{\db^{\frac{d}{2} + 1}_{2,1}} 
        + \norm{u_0}^h_{\db^{\frac{d}{2}}_{2,1}}
        \\
        & + \sum_{k=1}^d
         \Bigc{\norm{ ( x_k a_0, x_k u_0 ) }^l_{\db^{\frac{d}{2}}_{2,1} }
        +
        \norm{x_k a_0}^h_{ \db^{\frac{d}{2} + 1}_{2,1}}
        + \norm{x_k u_0}^h_{\db^{\frac{d}{2}}_{2,1}} }
        +
        \norm{(a_0,u_0)}_{ \db^{-\frac{d}{2}}_{2,\infty} }
        \leq c,
    \end{align*}
    then \eqref{CNS} has a unique global-in-time solution $(a,u)$ in the space $S$ defined above. Also, there exists a constant $C = C(d,\mu, P, j_0)$ such that
    \begin{align*}
        \norm{(a,u)}_{S} \leq C S_{0}.
    \end{align*}
\end{prop}

In this paper, we will be making use of the $d=3$ case of the above two propositions.

\begin{rem}
    Note that the initial data in Proposition \ref{3rd paper main thm} also satisfies the conditions for Proposition \ref{dx2017}.
\end{rem}

Returning to problem \eqref{CNS},
by applying the orthogonal projections $\mathcal{P}$ and $\mathcal{Q}$ onto the divergence and curl-free fields, respectively, and setting $\al \coloneqq P'(1)$ and $\nu \coloneqq \lambda + 2\mu$, we get the system
\begin{align} \label{div and curl free split}
    \begin{cases}
    \pt_t a + \div (\mathcal{Q} u) = f & \text{ in } \R_{> 0} \times \R^3, \\
    \pt_t \Q u - \nu \Delta \Q u + \al \nabla a = \Q g & \text{ in } \R_{> 0} \times \R^3, \\
    \pt_t \P u - \mu \Delta \P u = \P g & \text{ in } \R_{> 0} \times \R^3.
    \end{cases}
\end{align}
We set 
\[
v \coloneqq |D|^{-1} \div ( u), \text{ where }
|D|^s u  \coloneqq \F^{-1} \Big{[} |\xi|^s \hat{u}\Big{]}, \ s\in\R.
\]
We note that one can obtain $v$ from $\Q u$ by a Fourier multiplier of homogeneous degree zero. Thus, bounding $v$ is equivalent to bounding $\Q u$ in any Besov space (see Proposition \ref{fm est for besov}).

We note that we can set $\al = \nu = 1,$ without loss of generality, since the following rescaling
\begin{align*}
    a(t,x) = \tilde{a}\Bigc{\frac{\al}{\nu} t, \frac{\sqrt{\al}}{\nu} x }, \quad u(t,x) = \sqrt{\al} \ \tilde{u} \Bigc{ \frac{\al}{\nu} t , \frac{\sqrt{\al}}{\nu} x }
\end{align*}
ensures that $(\tilde{a}, \tilde{u})$ solves \eqref{div and curl free split} with $\al = \nu = 1.$
Thus we get that $(a, v)$ solves the following system:
\begin{align} \label{lin system av}
    \begin{cases}
    \pt_t a + |D| v = f & \text{ in } \R_{> 0} \times \R^3, \\
    \pt_t v - \Delta v - |D| a = h \coloneqq |D|^{-1} \div (g) & \text{ in } \R_{> 0} \times \R^3.
    \end{cases}
\end{align}



In \cite{ohaodha-iwabuchi2023(LinCNS)}, the authors considered the homogeneous case, where $f = h = 0$, which gives the system
\begin{align}
\label{homogeneous case}
    \begin{cases} 
    \pt_t a + |D| v = 0 & \text{ in } (0,\infty) \times \R^3, \\
    \pt_t v - \Delta v - |D| a = 0 & \text{ in } (0,\infty) \times \R^3.
    \end{cases}
\end{align}
Taking the Fourier transform over space $x$, we can write the above system as
\begin{align} \label{ODE}
    \frac{\dd }{\dd t}
    \begin{bmatrix}
    \ha \\
    \hv
    \end{bmatrix}
    =
    M_{|\x|} \begin{bmatrix}
    \ha \\
    \hv
    \end{bmatrix}
    ,
    \quad \text{with}
    \quad M_{|\x|} \coloneqq 
    \begin{bmatrix}
    0 & -|\x| \\
    |\x| & -|\x|^2
    \end{bmatrix}.
\end{align}
Then we may write the following formula for the solution to \eqref{homogeneous case}:
\begin{align*}
    \begin{bmatrix}
    a (t) \\
    v (t)
    \end{bmatrix} 
    & = 
    e^{t M(D) } 
    \begin{bmatrix}
    a_0 \\
    v_0
    \end{bmatrix}
    =
    \F^{-1}
    \Big{[}
    e^{t M_{|\x|}} 
    \begin{bmatrix}
    \ha_0 \\
    \hv_0
    \end{bmatrix}
    \Big{]}
    .
\end{align*}


The authors obtained in \cite{ohaodha-iwabuchi2023(LinCNS)} the following sharp decay result for the linear solution:


\begin{prop} {\rm(}\cite{ohaodha-iwabuchi2023(LinCNS)}{\rm)} \label{2nd paper main thm}
Let $s \in \R$, $p \in [2,\infty],$ $q \in [1,\infty]$, and $t > 1$.
If $a_0,v_0 \in \db^{s}_{1,q} \cap \db^s_{p,q}$,
then there exists $C>0$ such that
\begin{align} \label{thm1 l1 linf}
        \bignorm{ 
        e^{t M(D) } 
        \begin{bmatrix}
        a_0 \\
        v_0
        \end{bmatrix}
        }_{\db^s_{p,q}}
        \leq 
        C t^{-\frac{3}{2} (1-\frac{1}{p}) - \frac{1}{2} (1 - \frac{2}{p}) }
        \bignorm{ 
        \begin{bmatrix}
        a_0 \\
        v_0
        \end{bmatrix}  }^l_{\dot{B}^{s}_{1,q}}
        +
        C e^{-t}
        \bignorm{  
        \begin{bmatrix}
        a_0 \\
        v_0
        \end{bmatrix}  }^h_{
        \dot{B}^{s}_{p,q} }.
\end{align}
Also, there exist $a_0, v_0$ such that, for all sufficiently large $t$,
\begin{align} \label{thm 1 linf below}
    \bignorm{ 
    e^{t M(D) } 
        \begin{bmatrix}
        a_0 \\
        v_0
        \end{bmatrix}
    }_{\infty} \geq C t^{-2}.
\end{align}
\end{prop}


\bigskip

Returning to the inhomogeneous case, we use Duhamel's principle to obtain the integral formula for the solution to \eqref{lin system av}:
\begin{align}
\label{integral formula for full soln}
    \begin{bmatrix}
    a (t) \\
    v (t)
    \end{bmatrix} 
    & = 
    e^{t M(D) } 
    \begin{bmatrix}
    a_0 \\
    v_0
    \end{bmatrix}
    +
    \int^t_0
    e^{(t-s) M(D) } 
    \begin{bmatrix}
    f(s) \\
    h(s)
    \end{bmatrix}
    \dd s.
\end{align}
We now give our main result for this paper. 
We extend the results of \cite{ohaodha-iwabuchi2023(LinCNS)} to the radially symmetric inhomogeneous case by proving that the solution to \eqref{integral formula for full soln} (and thus the solution to \eqref{CNS}) decays at the same rate as the low-frequency decay in \eqref{thm1 l1 linf}. 
We also show that the nonlinear term decays at a faster rate, equivalent to the decay of the first derivative of the solution proven in Proposition \ref{Hoff-Zumbrun}.

We let $E$ and $S_0$ be the norms of the initial data in Propositions \ref{Hoff-Zumbrun} and \ref{3rd paper main thm}, respectively.
We recall the definitions here: 
\begin{align*}
    E  \coloneqq &
    \bignorm{  
    \begin{bmatrix}
        a_0 \\
        m_0
    \end{bmatrix}
    }_1
    +
    \bignorm{
    \begin{bmatrix}
        a_0 \\
        m_0
    \end{bmatrix}}_{H^{ 1 + l }},
    \\
    S_{0}  \coloneqq &  
        \norm{(a_0,u_0)}^l_{\dot{B}^{\frac{d}{2} - 1}_{2,1}} 
        + \norm{a_0}^h_{\db^{\frac{d}{2} + 1}_{2,1}} 
        + \norm{u_0}^h_{\db^{\frac{d}{2}}_{2,1}}
        \\
        & + \sum_{k=1}^d
         \Bigc{\norm{ ( x_k a_0, x_k u_0 ) }^l_{\db^{\frac{d}{2}}_{2,1} }
        +
        \norm{x_k a_0}^h_{ \db^{\frac{d}{2} + 1}_{2,1}}
        + \norm{x_k u_0}^h_{\db^{\frac{d}{2}}_{2,1}} }
        +
        \norm{(a_0,u_0)}_{ \db^{-\frac{d}{2}}_{2,\infty} }
        .
\end{align*}

\begin{thm} \label{4th paper main thm}
Let $p \in [2,\infty].$
Let $(a_0,u_0)$ be initial data for \eqref{CNS} satisfying the conditions for Proposition \ref{3rd paper main thm} and for Proposition \ref{Hoff-Zumbrun} with $l = 9$.
I.e., we take $E + S_0 \ll 1$.
Also assume that
\[ 
    \norm{ |D|^{-1} \div (u_0) }_1 
    + 
    \bignorm{ |\cdot| \begin{bmatrix}
        a_0 \\ |D|^{-1} \div (u_0) 
    \end{bmatrix} }_1
    < \infty, 
\]
and that
for all $x\in\R^3$,
\[
a_0(x) = a_0(|x|), \quad u_0(x) = U_0 (|x|) \frac{x}{|x|}, 
\]
where $U_0 : [0,\infty) \to \R$.
Then there exists a constant $C = C (a_0, u_0 ) > 0$ such that, for all $t \geq 1$,
\begin{align} \label{p4 main thm full est}
        \bignorm{ 
        \begin{bmatrix}
        a(t) \\
        u(t)
        \end{bmatrix}  }_{p }
        \leq 
        C t^{-\frac{3}{2} (1-\frac{1}{p}) - \frac{1}{2} (1 - \frac{2}{p}) }
        .
\end{align}
Also,
\begin{align} \label{p4 main thm nonlin est}
    \bignorm{
    \int^t_0
    e^{(t-s) M(D) } 
    \begin{bmatrix}
    f(s) \\
    h(s)
    \end{bmatrix}
    \dd s
    }_{p }
    \leq
    C
    t^{-\frac{3}{2} (1-\frac{1}{p}) - \frac{1}{2} (1 - \frac{2}{p}) -\frac{1}{2}}
        .
\end{align}
Finally, there exist initial data $(a_0,u_0)$ and a constant $C_0 > 0$ such that, for all $t>1$ sufficiently large,
\begin{align} \label{p4 main thm bnd from below}
    \bignorm{ 
        \begin{bmatrix}
        a(t) \\
        u(t)
        \end{bmatrix}  }_{ \infty }
        \geq 
        C_0 t^{ -2 }
        .
\end{align}
\end{thm}

\bigskip

\begin{rem}
    The condition $l=9$ is taken so that we may exploit the decay rates in Proposition \ref{Hoff-Zumbrun} up to the third derivative, which makes our argument simpler.
\end{rem}

The above theorem can in fact be refined by obtaining the same bound from above for the norm of the Besov space $\db^0_{p,1}$, which is stronger than the Lebesgue space $L^p$ (see Proposition \ref{Besov embedding}).
The bound from below can also be obtained for the $\db^0_{ \infty , \infty }$-norm, which is smaller than the $L^\infty$-norm.
In fact, our proof of \eqref{p4 main thm nonlin est} for $p \neq 2$ is reliant on estimates of the Besov norm.
We give this result as a separate theorem below.



\begin{thm} \label{4th paper 2nd main thm Besov}
Let $p \in [2,\infty],$ and $t > 0$.
Let $(a_0,u_0)$ be initial data for \eqref{CNS} satisfying the conditions for Theorem \ref{4th paper main thm}.
Then there exists a constant $C = C (a_0, u_0) > 0$ such that, for all $t \geq 1$,
\begin{align} \label{p4 2nd main thm full est}
        \bignorm{ 
        \begin{bmatrix}
        a(t) \\
        u(t)
        \end{bmatrix}  }_{ \db^{0}_{p,1} }
        \leq 
        C t^{-\frac{3}{2} (1-\frac{1}{p}) - \frac{1}{2} (1 - \frac{2}{p}) }
        .
\end{align}
Also,
\begin{align} \label{p4 2nd main thm nonlin est}
    \bignorm{
    \int^t_0
    e^{(t-s) M(D) } 
    \begin{bmatrix}
    f(s) \\
    h(s)
    \end{bmatrix}
    \dd s
    }_{ \db^{0}_{p,1} }
    \leq
    C
    t^{-\frac{3}{2} (1-\frac{1}{p}) - \frac{1}{2} (1 - \frac{2}{p}) -\frac{1}{2}}
        .
\end{align}
Finally, there exist initial data $(a_0,u_0)$ and a constant $C_0>0$ such that, for all $t>1$ sufficiently large,
\begin{align} \label{p4 2nd main thm bnd from below}
    \bignorm{ 
        \begin{bmatrix}
        a(t) \\
        u(t)
        \end{bmatrix}  }_{ \db^{0}_{\infty, \infty } }
        \geq 
        C_0 t^{ -2 }
        .
\end{align}
\end{thm}

\noindent {\bf Notation}

We obtain the following eigenvalues for $M_{|\x|} \coloneqq 
    \begin{bmatrix}
    0 & -|\x| \\
    |\x| & -|\x|^2
    \end{bmatrix} $, which differ between high and low frequencies:
\begin{align*}\lambda_{\pm}(\x) \coloneqq 
    \begin{cases}
    -\frac{|\x|^2}{2} \Bigc{1 \pm i { 
    \sqrt{\frac{4}{|\x|^2} - 1 } } }, & \text{ for } |\x| <2, \\
    -\frac{|\x|^2}{2} \Bigc{1 \pm { 
    \sqrt{ 1 - \frac{4}{|\x|^2} } } }, & \text{ for } |\x| >2.
    \end{cases}
\end{align*}

Throughout this paper, we will also use the following notation for the semigroup $e^{t \lambda_\pm(D)}$: we define the function $\G_\pm : [0,\infty) \times \R^3 \to \R$ such that, for all $f \in \S'_h$,
\[
\G_\pm(t)\ast f = e^{t \lambda_\pm (D)} f
=
\F^{-1} \Bigf{
e^{ t \lambda_{\pm} (\x) } \hat{f}
}
.
\]
Similarly we let $\G : [0,\infty) \times \R^{3} \to \R^{ 2 \times 2 }$ denote the function such that for all $f \in (\S'_h)^3$,
\[
\G (t) \ast f = e^{t M (D)} f 
=
\F^{-1} \Bigf{
e^{ t M_{|\x|} } \hat{f}
}
.
\]




\section{Preliminaries}

In the following section, we write several lemmas and definitions only in the 3-dimensional case.

\begin{definition}(The Fourier Transform)
For a function, $f$, we define the Fourier transform of $f$ as follows:
\begin{align*}
\F[f](\x) \coloneqq \hat{f}(\x) \coloneqq \frac{1}{(2\pi)^{3/2}}\intr e^{-i x\cdot\x} f(x) \dd x.
\end{align*}
The inverse Fourier transform is then defined as 
\begin{align*}
\F^{-1}[\hat{f}](x) \coloneqq \frac{1}{(2\pi)^{3/2}} \intr e^{i x \cdot \xi} \hat{f}(\x) \dd \x.
\end{align*}
For the purpose of calculating inequalities, we will frequently omit the factor of $1/(2\pi)^{3/2}$.
\end{definition}

\begin{definition}(Orthogonal Projections on the divergence and curl-free fields)
The projection mapping $\P$ is a matrix with each component defined as follows for $i,j \in \{1, 2, 3\}$:
\begin{align*}
(\P)_{i,j} \coloneqq \delta_{i,j} + (-\Delta)^{-1} \pt_i \pt_j.
\end{align*}
We then define $\Q \coloneqq 1 - \P$.
For $f \in (\dot{B}^s_{p,q}(\R^3))^3$, with $s \in \R$, and $p, q \in [1,\infty]$, we may write
\begin{align*}
\P f \coloneqq (1 + (-\Delta)^{-1} \nabla \div) f.
\end{align*}
\end{definition}

We next write some key properties of Besov spaces, whose proofs can be found in \cite{danchinbook}.

\begin{prop} \label{Besov embedding}
    Let $p\in[1,\infty]$.
    Then we have the following continuous embeddings:
    \begin{align*}
        \db^0_{p,1} \hookrightarrow L^p \hookrightarrow \db^0_{p,\infty}.
    \end{align*}
\end{prop}

\begin{prop}
    Let $s\in\R,$ $1\leq p_1 \leq p_2 \leq \infty$, and $1\leq r_1 \leq r_2 \leq \infty$.
    Then
    \begin{align*}
        \db^s_{p_1,r_1} \hookrightarrow \db^{s - 3 ( \frac{1}{p_1} - \frac{1}{p_2} )}_{p_2,r_2}.
    \end{align*}
\end{prop}

\begin{prop}
    Let $1\leq p \leq q \leq \infty$.
    Then
    \begin{align*}
        \db^{ \frac{3}{p} - \frac{3}{q} }_{p,1} \hookrightarrow L^q.
    \end{align*}
    Also, if $p < \infty$, then $\db^{ \frac{3}{p} }_{p,1}$ is continuously embedded the space $C_0$ of bounded continuous functions vanishing at infinity.
\end{prop}

For the next proposition, we introduce the notation $F(D)u \coloneqq \F^{-1} [ F (\cdot) \hu (\cdot) ]$.
\begin{prop}\label{fm est for besov} {\rm (}Fourier Multiplier Estimate{\rm )}
    Let $F$ be a smooth homogeneous function of degree $m$ on $\R^d \backslash \{0\}$ such that $F(D)$ maps $\mathcal{S}'_h$ to itself.
    Then
    \begin{align*}
        F(D) : \db^s_{p,r} \to \db^{s-m}_{p,r}.
    \end{align*}
    In particular, the gradient operator maps $\db^s_{p,r}$ to $\db^{s-1}_{p,r}$.
\end{prop}

\begin{prop} {\rm (}Composition Estimate{\rm )}
    Let $F : \R \to \R$ be smooth with $F(0) = 0$.
    Let $s>0$ and $1\leq p, r \leq \infty$.
    Then $F(u) \in \db^s_{p,r} \cap L^\infty$ for $u \in \db^s_{p,r} \cap L^\infty$, and there exists a constant $C = C( \norm{u}_{L^\infty}, F', s,p )>0$ such that
    \begin{align*}
        \norm{F(u)}_{ \db^s_{p,r} }
        \leq C
        \norm{ u }_{ \db^s_{p,r} }.
    \end{align*}
\end{prop}

\begin{prop} \label{p4 Danchin Besov-Linf Corollary}
    Let $u,v \in L^\infty \cap \db^s_{p,r}$, with $s>0$ and $1\leq p, r \leq \infty$.
    Then there exists a constant $C = C(p,s)>0$ such that
    \begin{align*}
        \norm{uv}_{ \db^s_{p,r} }
        \leq
        C \biggc{
        \norm{ u }_{L^\infty} \norm{ v }_{ \db^s_{p,r} }
        +
        \norm{ v }_{L^\infty} \norm{ u }_{ \db^s_{p,r} }
        }.
    \end{align*}
\end{prop}

We next discuss the existence of radial solutions.
\begin{prop} \label{p4 rad sym soln exist}
    Let $(a_0,u_0)$ satisfy the conditions of Theorem \ref{4th paper main thm}.
    Then the unique solution of \eqref{CNS} is radial.
    That is, for all $t>0$, and all $x \in \R^3$,
    \begin{align*}
        a(t,x) = a ( t , |x| ), \quad u(t,x) = U(t,|x|) \frac{x}{|x|},
    \end{align*}
    where $U : (0,\infty) \times [0,\infty) \to \R$.
\end{prop}

\begin{pf}
Let $R$ be a rotation matrix. 
We define
\begin{align*}
    a_R \coloneqq a(Rx), \quad u_R \coloneqq R^{-1} u (Rx)
\end{align*}
and observe that a radial solution is a solution that satisfies $(a_R,u_R) = (a,u)$.
We will prove that, if $(a,u)$ is a unique solution to \eqref{CNS} which is sufficiently regular (such that derivatives can be taken in the classical sense) with initial data $(a_0,u_0)$, then $(a_R , u_R)$ is a unique solution for \eqref{CNS} with initial data $(a_{0,R},u_{0,R})$.
Sufficient regularity of solutions is guaranteed by our setting $l=9$ when applying Theorem \ref{Hoff-Zumbrun}.

We note the following identities:
\begin{align*}
    (1) \ & \N (a_R) 
    = R^{-1} (\N a)(Rx),
    \\
    (2) \ & \div(u_R) 
    = (\div(u))(Rx).
    \\
    (3) \ & \N \Bigc{ (\div(u))(Rx) } 
    = R^{-1} (\N \div(u)) (Rx)
    \\
    (4) \ & (u_R \cdot \N) u_R
    = R^{-1} \Bigc{
    \bigc{(u \cdot \N) u } (Rx)
    }
    \\
    (5) \ & \Delta (u_R) 
    = 
    R^{-1} \bigc{\Delta u} (Rx).
\end{align*}
Using these, we see that if we apply the change of variables $x \to Rx$ to \eqref{CNS} and multiply the momentum equation by $R^{-1}$, then the equations for $(a,u)$ with initial data $(a_0,u_0)$ becomes the same equations for $(a_R,u_R)$ with initial data $(a_{0,R},u_{0,R})$.

Finally, suppose the initial data is radial as in Theorem \ref{4th paper main thm}, that is, for all $x\in\R^3$,
\[
a_0(x) = a_0(|x|), \quad u_0(x) = U_0 (|x|) \frac{x}{|x|}, 
\]
where $U : [0,\infty) \to \R$.
Then $(a_{0,R},u_{0,R}) = (a_0,u_0)$.
Thus, by the uniqueness of solutions in Theorem 
\ref{3rd paper main thm}, we get that
$(a_R,u_R) = (a,u)$, and thus the solution is radial.
\end{pf}

\bigskip

We will also make use of the following lemma:
\begin{lem} {\rm (\cite{danchinbook})} \label{p4 Danchin book Cor 1.39}
    Let $p\in(1,2]$. Then $L^p$ is continuously embedded in $\dot{H}^s$, with $s= \frac{3}{2} - \frac{3}{p}$.
\end{lem}

For the proof of the estimate \eqref{p4 main thm nonlin est} of the nonlinear term of the solution $(a,u)$, we will require a time-decay estimate of $(a,v)$ in the weighted $L^\infty$-norm.
We briefly explain our notation for weighted norms.
We write for a function $f$, and for $p\in[1,\infty]$,
\begin{align*}
    \norm{ x f }_p \coloneqq 
    \Bigc{ 
    \intr 
    | x f(x) |^p
    \dd x
    }^{\frac{1}{p} },
\end{align*}
and the meanings of the norms $\norm{|x|f}_p$ and $\norm{ x_k f }_p$ are similar.
We will also denote
\begin{align*}
    (xf) \ast g \coloneqq 
    \intr 
    y f(y) g(x-y)
    \dd y.
\end{align*}

Our weighted estimate requires the following lemma:
\begin{lem}\label{p4 rad sym ft ineq}
    Let $k \in \{ 1, 2, 3 \}$.
    Let $ f : \R^3 \to \R $ such that $f \in \mathcal{S}_h'$, $ f (x) = f (|x|)$ for all $x\in \R^3$, and $|\x| \hat{ f } \in L^1$.
    Then 
    \begin{align*}
        \norm{ x_k f }_\infty
        \leq 
        4\pi \int^\infty_0
        | \hat{ f } (\rho) | \rho 
        \dd \rho
    \end{align*}
\end{lem}

\begin{pf}
Without loss of generality, assume $k=3$.
Writing out the norm, we see that
\begin{align*}
    & \norm{ x_3 f }_\infty
    =
    \sup_{x\in\R^3}
    \Big{|} x_3 f (x) \Big{|} 
    \\
    & = 
    \sup_{x_3 \in \R }
    \Big{|} x_3 f ( x_3 e_3 ) \Big{|} 
\end{align*}
by the radial symmetry of $f$, where $e_3$ denotes the unit vector along the $x_3$-axis.
Next, we rewrite $f = \F^{-1} [\hat{ f }]$ and write out the inverse Fourier transform:
\begin{align*}
    & \F^{-1} [\hat{ f }] ( x_3 e_3 )
    =
    \int_{\R^3}
    e^{ i x_3 e_3 \cdot \x }
    \hat{ f } ( \x )
    \dd \x.
\end{align*}
We consider the dot product $ x_3 e_3 \cdot \x = |x_3| |\x| \cos{\t} $, where $\t$ is the angle between $x_3 e_3$ and $\x$, and thus also the angle between $\x$ and the $\x_3$-axis.
Thus, converting the integral coordinates to spherical coordinates $\x = (\rho, \t, \phi)$, we get
\begin{align*}
    & \int_{\R^3}
    e^{ i x_3 e_3 \cdot \x }
    \hat{ f } ( \x )
    \dd \x
    =
    \int^\infty_0
    \int^\pi_0
    \int^{2\pi}_0
    e^{ i |x_3| \rho \cos\t }
    \hat{ f } (\rho)
    \rho^2
    \sin\t
    \dd \phi
    \dd \t
    \dd \rho
    \\
    & =
    - \frac{ 4\pi }{ |x_3| }
    \int^\infty_0
    \hat{ f } \rho \sin( |x_3| \rho )
    \dd \rho,
\end{align*}
and so
\begin{align*}
    & \sup_{x_3 \in \R }
    \Big{|} x_3 f ( x_3 e_3 ) \Big{|} 
    =
    \sup_{x_3 \in \R }
    \Big{|} 4\pi  
    \int^\infty_0
    \hat{ f }(\rho) \rho \sin( |x_3| \rho )
    \dd\rho 
    \Big{|} 
    \\
    & \leq 
    4\pi 
    \int^\infty_0
    | \hat{ f } (\rho) | \rho 
    \dd \rho ,
\end{align*}
completing the proof of the lemma.
\end{pf}

\bigskip

The above lemma allows us to prove the following weighted $L^\infty$-estimate for $(a,v)$:
\begin{prop} \label{p4 wtd Linf est}
    Let $(a_0,u_0)$ satisfy the conditions of Theorem \ref{4th paper main thm}.
    Let $(a,v)$ be the associated solution to \eqref{lin system av}.
    Then there exists a constant $C = C(a_0,u_0) >0$ such that, for all $t>0,$
    \begin{align*}
        \bignorm{
        |x|
        \begin{bmatrix}
            a \\
            v
        \end{bmatrix} (t)
        }_\infty 
        \leq 
        C (t+1)^{ 
        -\frac{3}{4}
        }
        .
    \end{align*}
\end{prop}

\begin{rem}
    For this inequality, roughly speaking, we use boundedness of $xa,  xu$ in $\db^{s}_{2,1}$ from Proposition \ref{3rd paper main thm} with the inequality
    \begin{align*}
        \norm{ x e^{ t \lambda_\pm (D) } f }_\infty \leq 
        C t^{ - \frac{3}{4} }
        \Bigc{
        \norm{f}_2 + \norm{xf}_2
        }.
    \end{align*}
    Faster decay should be provable if sufficient decay results for $a,  u$ in weighted Besov spaces are obtained.
    However, for the proof of the main results in the present paper, the above inequality is sufficient.
\end{rem}

\begin{rem}
    Boundedness for $t\in(0,1)$ follows from the fact that, for all $t>0$,
    \[
    \bignorm{
        |x|
        \begin{bmatrix}
            a \\
            v
        \end{bmatrix} (t)
        }_\infty 
        \leq 
        C
        \Bigc{
        \sum_{k=1}^3
        \bignorm{
        x_k
        \begin{bmatrix}
            a \\
            u
        \end{bmatrix} (t)
        }_{ \db^{\frac{3}{2} }_{2,1}  } 
        +
        \norm{ u(t) }_{ \db^{\frac{3}{2} - 1}_{2,1} }
        }
        \leq C,
    \]
    where finiteness at the end follows from Proposition \ref{3rd paper main thm}.
    For the proof that follows, we thus focus on the $t \geq 1$ case.
\end{rem}

\begin{pf}
First, let us write the integral formula for the solution:
\begin{align*}
        \begin{bmatrix}
            a \\
            v
        \end{bmatrix} (t)
        =
        e^{t M(D)}
        \begin{bmatrix}
            a_0 \\
            v_0
        \end{bmatrix} 
        +
        \int^t_0
        e^{ (t-s) M(D)}
        F(s)
        \dd s,
\end{align*}
where 
\[
F(t) \coloneqq 
        \begin{bmatrix}
            f \\
            h
        \end{bmatrix} (t).
\]
Thus, we may split the norm of the solution as follows:
\begin{align} \label{p4 wtd l-inf est lin nonlin split}
    \bignorm{
        |x|
        \begin{bmatrix}
            a \\
            v
        \end{bmatrix} (t)
    }_\infty 
    \leq
    \sum^{3}_{k=1} \Biggc{
    \bignorm{
    x_k 
    \Bigc{e^{t M(D)}
        \begin{bmatrix}
            a_0 \\
            v_0
        \end{bmatrix} }
    }_\infty
        +
        \int^t_0
        \bignorm{
        x_k \Bigc{
        e^{ (t-s) M(D)}
        F(s)
        }
        }_\infty
        \dd s
    }.
\end{align}
We proceed by obtaining estimates for the linear and nonlinear terms separately, starting with the linear term.
We split linear term further into two terms, one where the $x_k$ is acting on the kernel of the semigroup, and one where it acts on the initial data:
\begin{align} \label{p4 lin term weighted est}
    \bignorm{
    x_k 
    \Bigc{e^{t M(D)}
        \begin{bmatrix}
            a_0 \\
            v_0
        \end{bmatrix} }
    }_\infty
    \leq 
    \bignorm{
    \Bigc{ x_k 
    \G(t)
    } \ast 
        \begin{bmatrix}
            a_0 \\
            v_0
        \end{bmatrix} 
    }_\infty
    +
    \bignorm{
    \G(t) \ast 
    \Bigc{
    x_k
        \begin{bmatrix}
            a_0 \\
            v_0
        \end{bmatrix} }
    }_\infty.
\end{align}
This is just a consequence of using the identity $x_k = x_k - y_k + y_k$ inside the integral in the formula for a convolution. I.e., for two functions $m,n$,
\begin{align*} 
& x_k (m \ast n) (x)
=
x_k
\intr
m(x-y) n(y)
\dd y
\\
& =
\intr
(x_k - y_k) m(x-y) n(y)
\dd y
+
\intr
m(x-y) y_k n(y)
\dd y
\\
& =
\bigc{( x_k m)\ast n} (x)
+
\bigc{ m \ast ( x_k n)} (x)
\end{align*}
for all $x\in\R^3$.
We consider the final term in \eqref{p4 lin term weighted est}. Looking at the Fourier transform, we see that
\begin{align*}
    & 
    \G(t) \ast 
    \Bigc{
    x_k
        \begin{bmatrix}
            a_0 \\
            v_0
        \end{bmatrix} }
    =
    \F^{-1} \Bigf{
    e^{t M( |\x| )}
    i \pt_k
    \begin{bmatrix}
            \ha_0 \\
            \hv_0
        \end{bmatrix} 
    }.
\end{align*}

We readily obtain by Proposition \ref{2nd paper main thm}
\begin{align*}
    & \bignorm{
    e^{t M(D)}
    \Bigc{
    x_k
        \begin{bmatrix}
            a_0 \\
            v_0
        \end{bmatrix} }
    }_\infty
    \leq C t^{-2} \Bigc{ 
    \bignorm{ x_k \begin{bmatrix}
            a_0 \\
            v_0
        \end{bmatrix}}_{1}
    +
    \bignorm{ x_k \begin{bmatrix}
            a_0 \\
            v_0
        \end{bmatrix}}^h_{ \db^0_{\infty, 1} }
    }
    \\
    & \leq 
    C t^{-2} \Bigc{ 
    \bignorm{ x_k \begin{bmatrix}
            a_0 \\
            v_0
        \end{bmatrix}}_{1}
    +
    S_0
    }
    \\
    & \leq 
    C t^{-2} 
    .
\end{align*}

We focus on the first term on the right-hand side of \eqref{p4 lin term weighted est}, namely 
$
\Bigc{ x_k 
\G(t)
} \ast 
    \begin{bmatrix}
        a_0 \\
        v_0
    \end{bmatrix} 
$.
Obtaining a bound for the $L^\infty$-norm here is more involved, as we cannot rely on a direct application of Proposition \ref{2nd paper main thm}. 
For brevity, we will just give a detailed proof of the decay result for the kernels of the individual semigroups $( x_k \G_\pm(t) ) \ast w$, where $w$ is a sufficiently regular generic function.
As discussed in \cite{ohaodha-iwabuchi2023(LinCNS)}, this is sufficient for frequencies away from $|\x| = 2$. 
See \cite{ohaodha-iwabuchi2023(LinCNS)} for discussion of how to estimate the kernel of the semigroup close to $|\x| = 2$.
The introduction of the weight $x_k$ does not change the strategy used there.
We start once again by splitting the semigroup kernel into its low, mid, and high frequencies:
\begin{align*}
    \norm{ ( x_k \G_\pm(t) ) \ast  w }_\infty
    \leq &
    \norm{ ( x_k \dot{S}_{-1} \G_\pm(t) ) \ast  w }_\infty
    +
    \norm{ ( x_k ( \dot{\Delta}_{0} + \dot{\Delta}_{1} + \dot{\Delta}_{2} )  \G_\pm (t) ) \ast w }_\infty
    \\
    &
    \ +
    \norm{ ( x_k ( 1 - \dot{S}_{2} ) \G_\pm (t) ) \ast w }_\infty.
\end{align*}
Let us start with the low-frequency estimate. 
We have
\begin{align*}
    \norm{ ( x_k \dot{S}_{-1} \G_\pm (t)  ) \ast w }_{\infty}
    \leq 
    \norm{ x_k \dot{S}_{-1} \G_\pm (t) }_\infty
    \norm{ w }_1.
\end{align*}
In order to estimate $\norm{ x_k \dot{S}_{-1} \G_\pm (t) }_\infty$, we apply Lemma \ref{p4 rad sym ft ineq}.
We get
\begin{align*}
    & \norm{ x_k \dot{S}_{-1} \G_\pm (t) }_\infty
    \leq C
    \int^2_0
    \Big{|}
    \hat{\psi}_{-1} (\rho)
    e^{ - t \frac{ \rho^2}{2} \Bigc{1 \pm i { 
    \sqrt{\frac{4}{\rho^2} - 1 } } } }
    \Big{|}
    \rho 
    \dd\rho 
    \\
    & \leq 
    C
    \int^\infty_0
    e^{ -t \frac{\rho^2}{2} }
    \rho 
    \dd \rho
    \\
    & =
    C t^{-1} \int^\infty_0
    e^{ -\frac{\rho^2}{2} } \rho
    \dd\rho
    \\
    & \leq 
    C t^{-1},
\end{align*}
where the last step is accomplished by a change of variables from $\rho$ to $t^{-\frac{1}{2}} \rho$.
The norm is also clearly bounded for small $t$, and thus we get
\begin{align*}
    & \norm{ x_k \dot{S}_{-1} \G_\pm (t) }_\infty
    \leq 
    C t^{-1},
\end{align*}
for all $t>0.$

We may bound the mid frequencies similarly, obtaining
\begin{align*}
    & \norm{ x_k 
    ( \dot{\Delta}_{0} + \dot{\Delta}_{1} + \dot{\Delta}_{2} ) 
    \G_\pm (t) }_\infty
    \\
    &
    \leq C
    \int^2_{1/2}
    e^{ -t \frac{\rho^2}{2} } \rho 
    \dd\rho
    + 
    \int^4_{2}
    e^{ -t \frac{\rho^2}{2} \Bigc{ 1 \pm 
    \sqrt{ 1-\frac{4}{\rho^2} } 
    } } \rho 
    \dd\rho
    \\
    & \leq 
    C e^{-ct},
\end{align*}
for some constant $c>0$ and all $t>0$.

Lastly, for the high frequencies, we use a different approach. 
To ensure boundedness of the less regular semigroup kernel $e^{t \lambda_-}$ at high frequencies, we apply an inverse Laplacian to it, with a compensatory Laplacian applied to $w$.
We then estimate the $L^\infty$-norm with two $L^2$-norms by Young's convolution inequality. We then get
\begin{align*}
    & \norm{ ( x_k ( 1 - \dot{S}_{2} )  \G_\pm (t)  ) \ast w }_\infty
    =
    \norm{ (-\Delta)^{-1} ( x_k ( 1 - \dot{S}_{2} )  \G_\pm (t)  )  \ast  (-\Delta) w }_\infty
    \\
    & \leq 
    C 
    \norm{ (-\Delta)^{-1} ( x_k ( 1 - \dot{S}_{2} ) \G_\pm(t) ) }_2
    \norm{ \Delta w }_2.
\end{align*}
Now that the semigroup kernel is in an $L^2$-norm, we may use the Plancherel therom to estimate it in Fourier space. 
We also make use of the following equality for the exponent of the semigroup kernel in Fourier space at high frequencies:
\begin{align} \label{p4 hi-freq exponent rewrite}
    {-t\frac{|\x|^2}{2} \bigc{ 1 \pm \sqrt{1 - {4}/{|\x|^2}}  }}
        =
        {-2t\bigc{1 \mp \sqrt{1 - {4}/{|\x|^2}}}^{-1}}
        =
        - t - \frac{4t}{|\x|^2} \bigc{1 \mp \sqrt{1 - 4 / |\x|^2}}^{-2}.
\end{align}
We thus may estimate
\begin{align*}
    & \norm{ (-\Delta)^{-1} ( x_k ( 1 - \dot{S}_{2} ) \G_\pm(t) ) }_2^2
    =
    \bignorm{  
    \frac{ \pt_k \Bigc{ ( 1 - \hat{\psi}_2 ) 
    e^{ t \lambda_\pm  }
    }}{|\x|^2}
    }_{L^2_\x}^2
    \\
    & \leq C
    \intr
    \Big{|}
    \frac{ \hat{\psi}_2' ( \x ) 
    e^{ -t \frac{|\x|^2}{2} \Bigc{
    1 \pm \sqrt{ 1-\frac{4}{|\x|^2} }
    } } 
    }{|\x|^2}
    \Big{|}^2
    \dd \x
    \\
    & \quad \quad \quad \quad 
    +
    C
    \intr
    \Bigg{|}
    \frac{ \Bigc{ 1 - \hat{\psi}_2 ( \x ) }
    e^{ -t \frac{|\x|^2}{2} \Bigc{
    1 \pm \sqrt{ 1-\frac{4}{|\x|^2} }
    } 
    }
    t \x_k \Bigc{
    (1 \pm \sqrt{ 1-\frac{4}{|\x|^2} } )
    \mp
    \frac{ 2 }{ |\x|^2  \sqrt{ 1-\frac{4}{|\x|^2} } }
    }
    }{|\x|^2}
    \Bigg{|}^2
    \dd \x
    \\
    & \leq C e^{-ct},
\end{align*}
for some $c>0$ and all $t>0$.
This completes the estimate for the linear term.

Returning to \eqref{p4 wtd l-inf est lin nonlin split}, we now consider the nonlinear term.
The estimate for the nonlinear term is complicated by the fact that
\begin{align*}
    & F \coloneqq 
    \begin{bmatrix}
        f \\
        h
    \end{bmatrix}
    =
    \begin{bmatrix}
        f \\
        |D|^{-1} \div(g)
    \end{bmatrix},
\end{align*}
and the presence of the $|D|^{-1} \div$, which may be thought of as a Riesz transform, is not easily ignored.
We must deal with the Riesz transform on $g$, while we can estimate the terms with $f$ similarly with fewer steps.

Thus, we shall focus on the estimate of 
\begin{align}
    & \notag
    \bignorm{ 
    x_k \Bigc{
    \G_\pm (t-s)
    \ast h (s)
    } 
    }_\infty
    \\
    & \notag
    \leq 
    \bignorm{ 
    \Bigc{ x_k 
    \dot{S}_{2}
    \G_\pm (t-s)
    } \ast 
    h (s)
    }_\infty
    +
    \bignorm{ 
    \Bigc{ x_k 
    ( 1 - \dot{S}_{2} )
    \G_\pm (t-s)
    } \ast 
    h (s)
    }_\infty
    \\
    & \ 
    +
    \bignorm{ 
    \dot{S}_{2} 
    \G_\pm (t-s) \ast 
    \Bigc{ x_k 
    h (s)
    }
    }_\infty
    +
    \bignorm{ 
    ( 1 - \dot{S}_{2} )
    \G_\pm (t-s) \ast 
    \Bigc{ x_k 
    h (s)
    }
    }_\infty.
    \label{p4 wtd est nonlin term breakdown}
\end{align}
Starting with the first term on the right-hand side, we estimate as follows.
First, let $\epsilon \in (0,1/4]$. Then,
\begin{align*}
    & \bignorm{ 
    \Bigc{ x_k 
    \dot{S}_{2}
    \G_\pm (t-s)
    } \ast 
    \dot{S}_3
    h (s)
    }_\infty
    =
    \bignorm{ 
    |D|^{ - 2\epsilon }
    \Bigc{ x_k 
    \dot{S}_{2}
    \G_\pm (t-s)
    } \ast
    |D|^{ 2\epsilon }
    \dot{S}_3
    h (s)
    }_\infty
    \\
    &\leq 
    \bignorm{ 
    |D|^{- 2\epsilon }
    \Bigc{ x_k 
    \dot{S}_{2}
    \G_\pm (t-s)
    }
    }_\infty
    \bignorm{
    |D|^{ 2\epsilon }
    \dot{S}_3
    h (s)
    }_1.
\end{align*}
Focusing on the left norm with the semigroup kernel first, we note that
\begin{align*}
    & |D|^{- 2\epsilon }
    \Bigc{ x_k 
    \dot{S}_{2}
    \G_\pm (t-s)
    }
    =
    \F^{-1}
    \Bigf{
    \frac{1}{|\x|^{ 2\epsilon }}
    \pt_k
    \Bigc{
    \hat{\psi}_{2}
    e^{(t-s) \lambda_\pm }
    }
    }
    \\
    & =
    \F^{-1}
    \Bigf{
    \pt_k
    \Bigc{
    \frac{1}{|\x|^{ 2\epsilon }}
    \hat{\psi}_{2}
    e^{(t-s) \lambda_\pm }
    }
    }
    -
    \F^{-1}
    \Bigf{
    \pt_k
    \Bigc{
    \frac{1}{|\x|^{ 2\epsilon }}
    }
    \hat{\psi}_{2}
    e^{(t-s) \lambda_\pm }
    }.
\end{align*}
Using Lemma \ref{p4 rad sym ft ineq} for the first term above, and directly estimating the second term, we get
\begin{align*}
    & \bignorm{
    |D|^{- 2\epsilon }
    \Bigc{ x_k 
    \dot{S}_{2}
    \G_\pm (t-s)
    }
    }_\infty
    \leq C
    \int^\infty_0
    \rho^{ 1 - 2\epsilon } 
    \hat{\psi}_{2} (\rho)
    e^{ -(t-s) \frac{\rho^2}{2} }
    \dd \rho
    \\
    & \leq
    C (t-s)^{ - 1 + \epsilon }.
\end{align*}
Next, for the norm of $|D|^{ 2\epsilon } \dot{S}_3 h (s),$ we write
\begin{align*}
    & \bignorm{
    |D|^{ 2\epsilon }
    \dot{S}_3
    h (s)
    }_1
    = 
    \bignorm{
    |D|^{ 2\epsilon }
    \dot{S}_3
    |D|^{-1} \div(g) (s)
    }_1
    \\
    & \leq 
    C 
    \norm{
    g (s)
    }_{\db^{ 2\epsilon }_{1,1}}^l
    \\
    & \leq
    C
    \norm{ g(s) }_1.
\end{align*}

We may next bound the second norm on the right-hand side of inequality \eqref{p4 wtd est nonlin term breakdown} in the same way that we estimated the mid and high-frequency parts of the linear term.
We have 
\begin{align*}
    & \bignorm{ 
    \Bigc{ x_k 
    ( 1 - \dot{S}_{2} )
    \G_\pm (t-s) 
    } \ast 
    h (s)
    }_\infty
    \\ &
    \leq C 
    \bignorm{ 
    \Delta^{-1}
    \Bigc{ 
    x_k 
    ( 1 - \dot{S}_{2} )
    \G_\pm (t-s)
    }
    }_2
    \bignorm{
    \Delta
    |D|^{-1} \div(g) (s)
    }_2
    \\
    & \leq C e^{-c (t-s) } \norm{ \Delta g(s) }_2.
\end{align*}

Moving onto the third norm on the right-hand side of inequality \eqref{p4 wtd est nonlin term breakdown}, 
\begin{align*}
    & \bignorm{ 
    \dot{S}_2
    \G_\pm (t-s) \ast 
    \Bigc{ x_k 
    h (s)
    }
    }_\infty
    \leq C
    \bignorm{ 
    \dot{S}_2
    \G_\pm (t-s)
    }_2
    \bignorm{
    \Bigc{ x_k 
    h (s)
    }
    }_2
    \\
    &
    \leq C (t-s)^{-\frac{3}{4}}
    \Bigc{
    \norm{ |D|^{-1} g (s) }_2
    +
    \norm{ x_k g(s) }_2
    }
    \\
    & \leq C (t-s)^{-\frac{3}{4}}
    \Bigc{
    \norm{  g (s) }_{\frac{6}{5}}
    +
    \norm{ x_k g(s) }_2
    },
\end{align*}
where the final step is obtained using Lemma \ref{p4 Danchin book Cor 1.39}.

Finally, for the last norm in \eqref{p4 wtd est nonlin term breakdown}, we simply apply the high-frequency estimate of Proposition \ref{2nd paper main thm}.
\begin{align*}
    & \bignorm{ 
    ( 1 - \dot{S}_2 )
    \G_\pm (t-s) \ast 
    \Bigc{ x_k 
    h (s)
    }
    }_\infty
    \leq C
    e^{-c(t-s)}
    \norm{ x_k g(s) }^h_{\db^0_{\infty,1}}.
\end{align*}
Combining the above inequalities for the linear term and nonlinear terms with $g$, and similar inequalities for the terms with $f$, we thus arrive at
\begin{align}
\notag
    & \bignorm{
        |x|
        \begin{bmatrix}
            a \\
            v
        \end{bmatrix} (t)
    }_\infty 
    \leq
    \sum^{3}_{k=1} \Biggc{
    \bignorm{
    x_k 
    \Bigc{e^{t M(D)}
        \begin{bmatrix}
            a_0 \\
            v_0
        \end{bmatrix} }
    }_\infty
        +
        \int^t_0
        \bignorm{
        x_k \Bigc{
        e^{ (t-s) M(D)}
        F(s)
        }
        }_\infty
        \dd s
    }
    \\
    \notag
    &
    \leq C t^{ 
    -1
    }
    +
    \int^t_0
    (t-s)^{-1} \norm{ f(s) }_1
    +
    (t-s)^{ - 1 + 2 \epsilon } \norm{ g(s) }_{ 1 }
    \dd s
    \\
    \notag
    & +
    \int^t_0
    e^{ -c(t-s) } 
    \Bigc{
    \norm{\Delta f(s)}_2
    +
    \norm{\Delta g(s)}_2
    }
    \dd s
    \\
    \notag
    & +
    \int^t_{0}
    (t-s)^{-\frac{3}{4}}
    \Bigc{
    \norm{ x_k f(s)}_2
    +
    \norm{ g(s) }_{ \frac{6}{5} }
    +
    \norm{ x_k g(s)}_2
    }
    \dd s
    \\
    \notag
    & +
    \int^t_0
    e^{-c(t-s)}
    \Bigc{
    \norm{ x_k f(s)}^h_{\db^0_{\infty,1}}
    +
    \norm{ x_k g(s)}^h_{\db^0_{\infty,1}}
    }
    \dd s
    \\
    & \leq C
    t^{-\frac{3}{4}}
    +
    \int^t_0
    e^{ - c(t-s) }
    \norm{u(s)}_{ \db^{\frac{3}{2} + 2 }_{2,1} }
    \bignorm{
        | x |
        \begin{bmatrix}
            a \\
            v
        \end{bmatrix} (s)
    }_\infty 
    \dd s,
    \label{p4 wtd est before gronwall step}
\end{align}
where the final step is obtained by splitting $f$ and $g$ with simple inequalities such as H\"older's inequality and Proposition \ref{p4 Danchin Besov-Linf Corollary} and then applying Propositions \ref{Hoff-Zumbrun}, \ref{dx2017}, and \ref{3rd paper main thm}.
The remaining integral term in \eqref{p4 wtd est before gronwall step} emerges from the following estimate of the last term in $g$:
\begin{align*}
    & \bignorm{
    \frac{ x_k a(s) }{1+a(s)}
    \mathcal{A} u(s)
    }_{\db^0_{\infty,1}}
    \\
    & \leq C
    \Bigc{
    \norm{ x_k a }_{\db^{\frac{3}{2}}_{2,1}}
    \norm{ \mathcal{A} u }_\infty
    +
    \norm{ x_k a }_\infty
    \norm{ \mathcal{A} u }_{\db^{\frac{3}{2}}_{2,1}}
    }
    \\
    & \leq C
    \Bigc{
    (s+1)^{ -\frac{5}{2} }
    +
    \norm{ x_k a }_\infty
    \norm{ \mathcal{A} u }_{\db^{\frac{3}{2}}_{2,1}}
    },
\end{align*}
Where $\norm{ x_k a (s) }_{\db^{\frac{3}{2}}_{2,1}}$ is bounded by Proposition \ref{3rd paper main thm} and $\norm{ \mathcal{A} u (s) }_\infty \leq C (s+1)^{ -\frac{5}{2} }$ by Proposition \ref{Hoff-Zumbrun}.
Finally, applying Gr\"onwall's inequality to \eqref{p4 wtd est before gronwall step} completes the proof.
\end{pf}




\section{Proof of Main Result}

We recall the problem under consideration:
\begin{align*}
    \begin{cases} 
    \pt_t a + |D| v = f & \text{ in } (0,\infty) \times \R^3, \\
    \pt_t v - \Delta v - |D| a = h & \text{ in } (0,\infty) \times \R^3.
    \end{cases}
\end{align*}

We write the nonlinear terms again for clarity:
\begin{align*}
    & f \coloneqq -\div(au),
    \\
    & h \coloneqq 
    |D|^{-1} \div \bigc{ -u\cdot\N u - \frac{a}{1+a} \mathcal{A} u - \beta(a) \N a } ,
\end{align*}
    where
\[ \beta(a) \coloneqq  \frac{P'(1+a)}{1+a} - P'(1).\]

We have the following integral formula for the solution:
\begin{align*}
    \begin{bmatrix}
    a (t) \\
    v (t)
    \end{bmatrix} 
    & = 
    e^{t M(D) } 
    \begin{bmatrix}
    a_0 \\
    v_0
    \end{bmatrix}
    +
    \int^t_0
    e^{(t-s) M(D) } 
    \begin{bmatrix}
    f(s) \\
    h(s)
    \end{bmatrix}
    \dd s.
\end{align*}

%
%
%
%

%
%
%
%
%
%
%
%

Proposition \ref{2nd paper main thm} gives us the estimate we need for the linear term. We now focus our attention on the nonlinear term.
\begin{prop} \label{p4 nonlin L2 est}
    Let $(a_0,u_0)$ satisfy the conditions of Theorem \ref{4th paper main thm}. There exists $C=C(Y)>0$ such that, for all $t >0,$
    \begin{align}
        \bignorm{ \NL }_{ 2 }
        \leq
        C t^{-\frac{3}{2} (1 - \frac{1}{2}) - \frac{1}{2}}
        .
    \end{align}
\end{prop}

\begin{pf}
For simplicity, we will only explicitly write the proof for the norm with the semigroup $e^{(t-s) \lambda_\pm(D) } $ instead of the whole matrix $e^{(t-s) M (D) } $.

We notice that, since 
\[ u(t,x) = U(t,|x|) \frac{x}{|x|}, \]
where $U : [0,\infty) \times [0,\infty) \to \R$ by Proposition \ref{p4 rad sym soln exist},
we may write
\[ u \cdot \N u = \N \Bigc{\frac{U^2}{2}}. \]
Also, using a Taylor expansion of $P$ and $1/(1+a)$, we can formally write $\beta(a) \N a$ as a series:
\[
\sum^\infty_{j=2} C_j \N(a^j),
\]
where, for sufficiently smooth $P,$ the sequence of constants $\{ C_j \}_{j\in\Nat}$ is bounded.

Thus, the terms $\div(au),$ $u\cdot \N u$, and $\beta(a) \N a$ may all be written in a `divergence form,' where a derivative operates on the whole term. This will allow us to `transfer' the derivative from the nonlinear term to the semigroup when taking estimates.
To show what we mean, we prove explicitly how the norm containing $\div(au)$ is estimated. 

We split the time interval into two halves and start with the `upper' half $\int^t_{t/2} ... \dd s$.
We will also estimate the low-frequencies and high-frequencies separately. 
Following \cite{ohaodha-iwabuchi2023(LinCNS)}, we choose $\dot{S}_2$ as our low-frequency cut-off function.
Let $t>0$. Then, starting with the low-frequency estimate, we get
\begin{align*}
    & \bignorm{ 
    \dot{S}_2 
    \int^t_{t/2}
    e^{-(t-s) \lambda_\pm (D)} \div(a(s) u(s))
    \dd s}_2
    \leq
    C
    \sum_{k=1}^3
    \int^t_{t/2}
    \norm{ \dot{S}_2 \G_\pm (t-s) }_2
    \norm{\pt_k (a(s) u_k(s))}_1
    \dd s
    \\
    & \leq C
    \int^t_{t/2}
    \norm{ \dot{S}_2 \G_\pm (t-s) }_2
    \bigc{
    \norm{\pt_k a}_2 \norm{u_k}_2
    +
    \norm{a}_2 \norm{\pt_k u_k}_2
    }
    \dd s
    \\
    & \leq C
    \int^t_{t/s}
    (t-s)^{ - \frac{3}{4}} E^2 (s+1)^{ -2 }
    \dd s
    \\
    & \leq C
    E^2 t^{-\frac{7}{4}},
\end{align*}
where the last step comes from a simple $L^2$-estimate of $\dot{S}_2 \G_\pm(t) $, which decays at the same rate of the heat kernel in $L^2$, and from applying Proposition \ref{Hoff-Zumbrun} to $(a,u)(s),$ recalling that
\[
    E \coloneqq 
    \bignorm{  
    \begin{bmatrix}
        a_0 \\
        m_0
    \end{bmatrix}
    }_1
    +
    \bignorm{
    \begin{bmatrix}
        a_0 \\
        m_0
    \end{bmatrix}}_{H^{ 1 + l }},
\]
where $m_0 \coloneqq \rho_0 u_0$, and we have taken $l=9$.
Next, we look at the high-frequency estimate.
Similarly to the high-frequency $L^p \textendash L^p$ estimate proven in \cite{ohaodha-iwabuchi2023(LinCNS)}, we observe that
\begin{align*}
        {-t\frac{|\x|^2}{2} \bigc{ 1 \pm \sqrt{1 - {4}/{|\x|^2}}  }}
        =
        {-2t\bigc{1 \mp \sqrt{1 - {4}/{|\x|^2}}}^{-1}}
        =
        - t - \frac{4t}{|\x|^2} \bigc{1 \mp \sqrt{1 - 4 / |\x|^2}}^{-2}.
    \end{align*}
We use this and Placherel's theorem to obtain
\begin{align*}
    & \bignorm{ 
    ( 1 - \dot{S}_2 ) 
    \int^t_{t/2}
    e^{-(t-s) \lambda_\pm (D)} \div(a(s) u(s))
    \dd s}_2
    \\
    &
    \leq 
    \sum_{k=1}^3
    \int^t_{t/2}
    \norm{  ( 1 - \hat{\psi}_2 ) e^{-(t-s) \lambda_\pm } \x_k (\ha(s) \hat{u}(s)) }_{L^2_\x}
    \dd s
    \\
    & \leq C
    \int^t_{t/2}
    e^{-(t-s)}
    \norm{ ( 1 - \hat{\psi}_2 ) e^{ - \frac{4(t-s)}{|\x|^2} \bigc{1 \mp \sqrt{1 - 4 / |\x|^2}}^{-2} } }_{L^\infty_\x}
    \norm{\N (a(s) u(s))}_2
    \dd s
    \\
    & \leq C
     E^2 t^{-\frac{ 11 }{4}}.
\end{align*}
Next, we look at the `lower' half of the time integral, $\int^{t/2}_0 ... \dd s$. In this case, in order to obtain our desired decay in the low-frequency estimate, we need to move the derivative on the nonlinear term across the convolution and onto the semigroup. 
\begin{align*}
    & \bignorm{ 
    \dot{S}_2 
    \int^{t/2}_0
    e^{-(t-s) \lambda_\pm (D)} \div(a(s) u(s))
    \dd s}_2
    \leq
    C
    \sum_{k=1}^3
    \int^{t/2}_0
    \norm{ \dot{S}_2 \pt_k \G_\pm (t-s) }_2
    \norm{ (a(s) u_k(s))}_1
    \dd s
    \\
    & \leq C \int^{t/2}_0
    (t-s)^{-\frac{3}{4} -\frac{1}{2} }
    E^2 (s+1)^{-\frac{3}{2}}
    \dd s
    \\
    & \leq C E^2 t^{ -\frac{3}{4} -\frac{1}{2}  }.
\end{align*}
Lastly, the high-frequency part decays so fast already that we obtain exponential decay by the same steps as on the upper half of the time integral.
\begin{align*}
    & \bignorm{ 
    ( 1 - \dot{S}_2 ) 
    \int^{t/2}_0
    e^{-(t-s) \lambda_\pm (D)} \div(a(s) u(s))
    \dd s}_2
    \\
    &
    \leq 
    \sum_{k=1}^3
    \int^{t/2}_0
    \norm{  ( 1 - \hat{\psi}_2 ) e^{-(t-s) \lambda_\pm } \x_k (\ha(s) \hat{u}(s)) }_{L^2_\x}
    \dd s
    \\
    & \leq C
    \int^{t/2}_0
    e^{-(t-s)}
    \norm{ ( 1 - \hat{\psi}_2 ) e^{ - \frac{4(t-s)}{|\x|^2} \bigc{1 \mp \sqrt{1 - 4 / |\x|^2}}^{-2} } }_{L^\infty_\x}
    \norm{\N (a(s) u(s))}_2
    \dd s
    \\
    & \leq C
    E^2 e^{- t/2 }  .
\end{align*}

The nonlinear terms containing $u\cdot \N u$ and $\beta(a) \N a$ are bounded by similar steps to the above. 
We thus move on to the final nonlinear term, 
\[
\int^{t}_0
    e^{-(t-s) \lambda_\pm (D)} |D|^{-1} \div( \frac{a}{1+a} \mathcal{A} u )
    \dd s,
\]
which presents a unique challenge, as it cannot be rewritten in a divergence form like the other nonlinear terms to transfer a derivative onto the semigroup. 
We thus need some other way of extracting the additonal $t^{-\frac{1}{2}}$ decay for the low-frequency estimate of the `lower' half of the time integral $\int^{t/2}_0 ... \dd s$.
All other estimates are similar to those we performed for nonlinear term with $\div(au)$, and so we focus on just this more difficult estimate.
Recall that we have set $ 2 \mu + \lambda = 1 $.
We note that, since $u$ is radial and thus curl-free, we may rewrite
\begin{align*}
& \mathcal{A} u = \mu \Delta u + (\lambda + \mu) \N \div(u) \\
& = 
\mu \Delta \Q u + (\lambda + \mu) \N \div( \Q u)
\\
& =
\mu \N \div( u) + (\lambda + \mu) \N \div(u)
\\
& =  |D| \N v.
\end{align*}
Thus, we are interested in the norm
\begin{align*}
    & \bignorm{
    \dot{S}_2 \int^{t/2}_{0}
    e^{ -(t-s)\lambda_\pm (D) }
    |D| \div \Bigc{
    \frac{a}{1+a} \N |D| v(s)
    }
    \dd s
    }_2
    \\ &
    =
    \bignorm{
    \dot{S}_2 \int^{t/2}_{0}
    \sum^{3}_{l=1}
    |D|^{-1} \pt_l
    e^{ -(t-s)\lambda_\pm (D) }
    \Bigc{
    \frac{a}{1+a} \pt_l |D| v(s)
    }
    \dd s
    }_2.
\end{align*}
Since $a$ and $v$ are both radially symmetric scalar functions, we get that the nonlinear term $\displaystyle \frac{a}{1+a} \pt_l |D| v(s)$ is antisymmetric, and thus its integral over space is zero. 
That is,
\begin{align*}
    \intr
    \frac{a(s, x)}{1+a(s, x)} \pt_l |D| v(s, x)
    \dd x
    = 
    0.
\end{align*}
We exploit this fact to place an extra derivative on the semigroup, in exchange for multiplication by the space variable (which we see in Proposition \ref{p4 wtd Linf est} behaves like an antiderivative) on the nonlinear term.
That is, we estimate 
\begin{align}
\notag
    & \bignorm{
    \dot{S}_2 \int^{t/2}_{0}
    |D|^{-1} \pt_l
    e^{ -(t-s)\lambda_\pm (D) }
    \Bigc{
    \frac{a}{1+a} \pt_l |D| v(s)
    }
    \dd s
    }_2
    \\
    \notag
    & = 
    \bignorm{
    \int^{t/2}_{0}
    \intr
    |D|^{-1} \pt_l
    ( \dot{S}_2 \G_\pm )(t-s, \cdot -y)
    \Bigc{
    \frac{a}{1+a} \pt_l |D| v(s,y)
    }
    \dd y
    \\ 
    \notag
    & 
    \quad \quad \quad \quad 
    -
    |D|^{-1} \pt_l
    ( \dot{S}_2 \G_\pm )(t-s, \cdot )
    \intr
    \frac{a}{1+a} \pt_l |D| v(s,y)
    \dd y
    \dd s
    }_2
    \\
    \notag
    & = 
    \bignorm{
    \int^{t/2}_{0}
    \intr
    |D|^{-1} \pt_l
    \biggc{
    ( \dot{S}_2 \G_\pm )(t-s, \cdot -y)
    - 
    ( \dot{S}_2 \G_\pm )(t-s, \cdot )
    }
    \Bigc{
    \frac{a}{1+a} \pt_l |D| v(s,y)
    }
    \dd y
    \dd s
    }_2
    \\
    \notag
    & =
    \bignorm{
    \int^{t/2}_{0}
    \intr
    \int^{1}_{0}
    |D|^{-1} \pt_l
    \N
    ( \dot{S}_2 \G_\pm )(t-s, \cdot - \theta y) 
    \cdot (-y)
    \dd \theta
    \Bigc{
    \frac{a}{1+a} \pt_l |D| v(s,y)
    }
    \dd y
    \dd s
    }_2
    \\
    \notag
    & \leq
    C t^{-\frac{3}{2}(1 - \frac{1}{2}) - \frac{1}{2}  }
    \int^{t/2}_{0}
    \bignorm{ | x | \frac{a}{1+a} \pt_l |D| v(s) }_{1}
    \dd s
    \\
    & \leq
    C t^{-\frac{3}{2}(1 - \frac{1}{2}) - \frac{1}{2}  }
    \int^{t/2}_{0}
    \bignorm{ | x | \frac{a}{1+a} }_{\infty}
    E (s+1)^{-\frac{1}{2}}
    \dd s, \label{p4 p=2 nonlin prop last ineq}
\end{align}
where the final step is an application of H\"older's inequality and Proposition \ref{Hoff-Zumbrun} for $p=1$.
By Proposition \ref{p4 wtd Linf est}, we know that
\[
\bignorm{ | x | \frac{a}{1+a} }_{\infty} \leq C(s+1)^{-\frac{ 3 }{ 4 }},
\]
for all $s>0$, and thus the time integral in \eqref{p4 p=2 nonlin prop last ineq} is bounded by a constant.
\end{pf}

\bigskip

The decay of the nonlinear term in the $L^\infty$-norm is not so easily proven, due to the Riesz transform in $h$.
We are forced in the end to bound the $L^\infty$-norm by the $\db^0_{\infty,1}$-norm.
We thus give the Besov-norm estimate next.

\begin{prop} \label{p4 nonlin Besov est}
    Let $(a_0,u_0)$ satisfy the conditions of Theorem \ref{4th paper main thm}. 
    There exists $C=C(Y)>0$ such that, for all $t >0,$
    \begin{align*}
        & \bignorm{ \NL }_{ \db^0_{2,1} }
        \leq
        C t^{-\frac{3}{2} (1 - \frac{1}{2}) - \frac{1}{2}}
        ,
        \\
        & \bignorm{ \NL }_{ \db^0_{\infty,1} }
        \leq
        C t^{ - 2 - \frac{1}{2}}
        .
    \end{align*}
\end{prop}

\begin{pf}
    The proof is similar to that of Proposition \ref{p4 nonlin L2 est}, but with individual dyadic blocks $\ddj$ replacing the low and high-frequency cut-offs $\dot{S}_2$ and $(1-\dot{S}_2)$. 
    We need only take care that the sum over $j$ is finite after obtaining our desired decay.
    Like the proof of Proposition \ref{p4 nonlin L2 est}, we will only explicitly consider the estimates of 
    $e^{t\lambda_\pm(D)}f$ and 
    $e^{t\lambda_\pm(D)} |D|^{ -1 } \div \Bigc{\frac{a}{1+a} \mathcal{A} u }$.
    We proceed in the same order as before, starting with the low frequencies in the upper half of the time integral $\int^{t}_{t/2} ... \dd s$.
    Let $t>0$.
    \begin{align*}
        & \sum_{j \leq 2}
        \bignorm{
        \ddj \int^t_{t/2}
        e^{(t-s) \lambda_\pm(D)}
        \div\Bigc{a(s)u(s)}
        \dd s
        }_2
        \leq 
        \sum_{j \leq 2} C
        \int^t_{t/2}
        \norm{ \ddj \G_\pm (t-s) }_2
        \norm{\N (a(s) u(s))}_1
        \dd s
        \\
        & \leq C 
        \int^t_{t/2}
        (t-s)^{ -\frac{3}{4} }  (s+1)^{-2}
        \dd s
        \\
        & \leq C t^{-2 + \frac{1}{4}},
    \end{align*} 
    which is more than fast enough.
    In the last step above, the norm of the semigroup kernel is estimated and the sum over $j$ taken as follows: 
    \begin{align}
        \notag
        & \sum_{ j \leq 2 }
        \norm{ \ddj \G_\pm (t-s) }_2
        =
        \sum_{ j \leq 2 }
        \Bigc{ \intr 
        \big{|} \hat{\phi} (2^{-j}\x) 
        e^{ -(t-s) \frac{|\x|^2}{2} \bigc{ 
        1 \pm i \sqrt{ \frac{4}{|\x|^2} - 1 }
        } } \big{|} ^2
        \dd \x }^{\frac{1}{2}}
        \\
        \notag
        & \leq 
         C 
        \sum_{ j \leq 2 }
        2^{\frac{3}{2}j} 
        \Bigc{ \intr 
        \big{|}  \hat{\phi} ( \x ) 
        e^{ -2^{2j} (t-s) \frac{|\x|^2}{2}  } \big{|} ^2
        \dd \x }^{\frac{1}{2}}
        \\
        \notag
        & \leq 
        C (t-s)^{-\frac{3}{4}} \sum_{ j \leq 2 }
        2^{\frac{3}{2}j} (t-s)^{\frac{3}{4}}
        e^{-c 2^{2j} t}
        \\
        & \leq 
        C (t-s)^{-\frac{3}{4}}.
        \label{p4 lo-freq sum over j fin}
    \end{align}
    The inequality for $p=\infty$ follows the exact same steps, except we use the following convolution inequality:
    \[
    \bignorm{
        \ddj \int^t_{t/2}
        e^{(t-s) \lambda_\pm(D)}
        \div\Bigc{a(s)u(s)}
        \dd s
        }_\infty
        \leq C 
        \int^t_{t/2}
        \norm{ \ddj \G_\pm (t-s) }_2
        \norm{\N (a(s) u(s))}_2
        \dd s.
    \]

    Next, we look at high frequencies $j>2$. 
    The steps are similar to the proof of the high-frequency $L^p \textendash L^p$ estimate in Proposition \ref{2nd paper main thm} (see \cite{ohaodha-iwabuchi2023(LinCNS)}).
    The semigroup kernel is rewritten in the same way as in \eqref{p4 hi-freq exponent rewrite}.
    We apply a Laplacian and inverted Laplacian.
    The inverted Laplacian ensures finiteness of the sum over $j>2$, and the Laplacian is readily absorbed by the nonlinear term.
    Let $p \in \{ 2, \infty \}$.
    \begin{align*}
        & \sum_{j>2}
        \bignorm{
        \ddj \int^t_{t/2}
        e^{(t-s) \lambda_\pm(D)} 
        (-\Delta)^{-1} (-\Delta)
        \div\Bigc{a(s)u(s)}
        \dd s
        }_p
        \\
        &
        \leq 
        \sum_{j>2} 2^{-2j} C \int^t_{t/2}
        e^{-(t-s)}
        \bignorm{
        \F^{-1}
        \Bigf{
        e^{-\frac{ 4 (t-s) }{ |\x|^2 }
        \bigc{ 1 \mp \sqrt{1 - \frac{4}{|\x|^2}} }^{-2}
        }
        \hp_j
        }
        }_1
        \norm{ \Delta \N (a(s) u(s) ) }_p
        \dd s
        \\
        & \leq C t^{ -\frac{3}{2} (1-\frac{1}{p}) -\frac{1}{2} (1 - \frac{2}{p}) - \frac{1}{2} }.
    \end{align*}
    The actual decay rate in the last step could be much faster, but we have bounded from above by our target decay for simplicity.
    The above $L^1$-norm is bounded by a constant as follows: for $j>2$,
    \begin{align*}
        \bignorm{
        \F^{-1} \Big{[}
        e^{ - \frac{4t}{|\x|^2} \bigc{ 1 \mp \sqrt{1 - \frac{4}{|\x|^2}} }^{-2} }
        \hp_j
        \Big{]}
        }_1
        & = 
        \bignorm{
        \F^{-1} \Big{[}
        e^{ - \frac{4t}{ 2^{2j} |\x|^2} \bigc{ 1 \mp \sqrt{1 - \frac{4}{ 2^{2j} |\x|^2}} }^{-2} }
        \hp_0
        \Big{]}
        }_1
        \\
        & \leq C
        \bignorm{
        e^{ - \frac{4t}{ 2^{2j} |\x|^2} \bigc{ 1 \mp \sqrt{1 - \frac{4}{ 2^{2j} |\x|^2}} }^{-2} }
        \hp_0
        }_{W^{2,2}}
        \\
        & \leq C. 
    \end{align*}
    
    Next, we consider the lower half of the time integral. 
    For low frequencies, we have 
    \begin{align*}
        & \sum_{j \leq 2}
        \bignorm{
        \ddj \int^{t/2}_0
        e^{(t-s) \lambda_\pm(D)}
        \div\Bigc{a(s)u(s)}
        \dd s
        }_p
        \leq
        \sum_{j \leq 2} C \int^{t/2}_0 
        \norm{ \ddj \N \G_\pm (t-s) }_p 
        \norm{a(s) u(s)}_1
        \dd s
        \\
        & \leq 
        C
        \int^{t/2}_0 
        (t-s)^{  -\frac{3}{2} (1-\frac{1}{p}) -\frac{1}{2} (1 - \frac{2}{p}) - \frac{1}{2}  }
        (s+1)^{- \frac{3}{2}}
        \dd s
        \\
        & \leq 
        C t^{ -\frac{3}{2} (1-\frac{1}{p}) -\frac{1}{2} (1 - \frac{2}{p}) - \frac{1}{2} },
    \end{align*}
    Where the sum over $j$ is handled similaryly to \eqref{p4 lo-freq sum over j fin} for the $p=2$ case.
    In the $p=\infty$ case, we estimate 
    by exploiting the heat-like component of the kernel of our semigroup as follows
    \[
    e^{ -t \frac{|\x|^2}{2} \bigc{
        1 \pm i \sqrt{ \frac{4}{|\x|^2} - 1 }
        } }
    =
    e^{ -t \frac{|\x|^2}{4} }
    e^{ -t \frac{|\x|^2}{4} \bigc{
        1 \pm 2i \sqrt{ \frac{4}{|\x|^2} - 1 }
        } }.
    \]
    We get
    \begin{align*}
        & \sum_{j\leq 2} \bignorm{
        \ddj \N \G_\pm (t-s)
        }_{\infty}
        \leq C \sum_{j\leq 2} 2^{ j }
        \bignorm{
        \ddj \F^{-1} \Bigf{ e^{-(t-s) \frac{|\x|^2}{4}} }
        }_{1}
        \bignorm{
        \dot{S}_3 \F^{-1} \Bigf{
        e^{ -(t-s) \frac{|\x|^2}{4} \bigc{
        1 \pm 2i \sqrt{ \frac{4}{|\x|^2} - 1 }
        } }
        }
        }_\infty
        \\
        & \leq C (t-s)^{ -2 -\frac{1}{2} }
        \sum_{j\leq 2}
        2^j (t-s)^{\frac{1}{2}} e^{ - c 2^{2j} (t-s) }
        \\
        & \leq C (t-s)^{ -2 -\frac{1}{2} }.
    \end{align*}
    The high-frequency part is estimated using the exact same steps as on the upper half of the time integral, but ending with exponential decay.

    Regarding the term with $e^{t\lambda_\pm(D)} |D|^{ -1 } \div \Bigc{\frac{a}{1+a} \mathcal{A} u }$, once again, the only part that is estimated differently from the terms in divergence form is the low-frequency part in the lower half of the time integral $\int^{t/2}_0 ... \dd s$.
    The same method as in Proposition \ref{p4 nonlin L2 est} is used to estimate the norm in this case, for both $p=2$ and $p=\infty$.
    \begin{align}
    \notag
    & 
    \sum_{ j \leq 2 }    
    \bignorm{
    \ddj \int^{t/2}_{0}
    |D|^{-1} \pt_l
    e^{ -(t-s)\lambda_\pm (D) }
    \Bigc{
    \frac{a}{1+a} \pt_l |D| v(s)
    }
    \dd s
    }_p
    \\
    \notag
    & =
    \sum_{ j \leq 2 }   
    \bignorm{
    \int^{t/2}_{0}
    \intr
    \int^{1}_{0}
    |D|^{-1} \pt_l
    \N
    ( \ddj \G_\pm )(t-s, \cdot - \theta y) 
    \cdot (-y)
    \dd \theta
    \Bigc{
    \frac{a}{1+a} \pt_l |D| v(s,y)
    }
    \dd y
    \dd s
    }_p
    \\
    \notag
    & \leq
    \sum_{ j \leq 2 }  
    \int^{t/2}_{0}
    \bignorm{
    |D|^{-1} \pt_l
    \N
    ( \ddj \G_\pm )(t-s) 
    }_p
    \bignorm{ |x|
    \frac{a}{1+a} \pt_l |D| v(s)
    }_1
    \dd s
    \\
    & \leq
    C t^{ -\frac{3}{2}(1 - \frac{1}{p}) - \frac{1}{2}(1-\frac{2}{p})  - \frac{1}{2}  }
    \int^{t/2}_{0}
    \bignorm{ | x | \frac{a}{1+a} \pt_l |D| v(s) }_{1}
    \dd s
    \label{1.7 pf last time int}
    \\
    \notag
    & \leq C t^{ -\frac{3}{2}(1 - \frac{1}{p}) - \frac{1}{2}(1-\frac{2}{p})  - \frac{1}{2}  },
    \end{align}
    where again, the sum over $j$ is taken by the same method as \eqref{p4 lo-freq sum over j fin}, and the time integral in \eqref{1.7 pf last time int} is bounded thanks to Proposition \ref{Hoff-Zumbrun} and Proposition \ref{p4 wtd Linf est}.
\end{pf}

\bigskip

%
%
%
%
%
%
%
%
%
%
%
%
%
%
%
%
%
%
%
%
%
%
%
%
%

%
%
%
%
%
%
%
%
%
%
%
%
%
%
%
%
%
%
%
%
%
%
%
%
%
%
%
%
%
%
%
%


\noindent{\bf Proof of Theorem \ref{4th paper main thm}.}
First, since $u$ is curl-free, we have for $v\coloneqq |D|^{-1} \div(u)$
\begin{align*}
    & u = \Q u
    = - (-\Delta)^{-1} \N \div(u)
    = - |D|^{-1} \N v,
\end{align*}
and thus, by Proposition \ref{fm est for besov}, we have that there exists a constant $C>0$ such that for all $t \geq 0$ and $p \in [1,\infty]$,
\[
C^{-1} \norm{ v(t) }_{\db^0_{p,1}}
\leq 
\norm{ u(t) }_{\db^0_{p,1}} 
\leq C
\norm{ v(t) }_{\db^0_{p,1}}.
\]
Estimates for $v$ in $\db^0_{p,1}$ thus imply estimates for $u$.
Estimates for $v$ in $L^2$ also clearly imply estimates for $u$ by Plancherel's theorem.

We obtain estimates of the linear term by Proposition \ref{2nd paper main thm}. 
Next, Proposition \ref{p4 nonlin L2 est} provides the nonlinear estimate \eqref{p4 main thm nonlin est} for $p=2$. The $p=\infty$ case is proven by Proposition \ref{p4 nonlin Besov est} and the fact that $\db^0_{\infty,1} \hookrightarrow L^\infty$. By interpolation, we obtain the nonlinear estimate for other values of $p\in (2,\infty)$.
Combining the estimates of the linear and nonlinear terms yields \eqref{p4 main thm full est}. 
Finally, similarly to the proof of optimality for Proposition \ref{2nd paper main thm} (see \cite{ohaodha-iwabuchi2023(LinCNS)}), 
if $a_0 = c e^{-|x|^2}$, for sufficiently small $c>0$,
then there exists $\delta > 0$ such that if $\norm{v_0}_{L^1 \cap \db^0_{\infty,1}} < \delta$, then 
\begin{align*}
    & \bignorm{ 
    \begin{bmatrix}
        a(t) \\ v(t)
    \end{bmatrix}
    }_\infty
    \geq 
    \bignorm{ 
    e^{t M(D)}
    \begin{bmatrix}
        a_0 \\ v_0
    \end{bmatrix}
    }_\infty
    -
    \bignorm{ 
    \int^t_0
    e^{(t-s) M(D)}
    \begin{bmatrix}
        f(s) \\ h(s)
    \end{bmatrix}
    \dd s
    }_\infty
    \\
    & \geq 
    C t^{-2} - C t^{-2 -\frac{1}{2}} 
    \bignorm{
    \begin{bmatrix}
        a_0 \\ u_0
    \end{bmatrix}
    }_Y
    \\
    & \geq C t^{-2},
\end{align*}
for all sufficiently large $t>0$.
\hfill\qed

\bigskip

\noindent{\bf Proof of Theorem \ref{4th paper 2nd main thm Besov}.}
We obtain estimates of the linear term by Proposition \ref{2nd paper main thm}. 
Next, Proposition \ref{p4 nonlin Besov est} provides the nonlinear estimate \eqref{p4 2nd main thm nonlin est}.
Combining the estimates of the linear and nonlinear terms yields \eqref{p4 2nd main thm full est}.
Finally, the bound from below follows similarly to that of Theorem \ref{4th paper main thm}, after applying Proposition \ref{p4 Appendix bnd from below b0infin}, proven in Appendix \ref{p4 appendix}.
\hfill\qed


\appendix
\section{ The Bound from Below }
\label{p4 appendix}

The bound from below for the linear term, as proven in~\cite{ohaodha-iwabuchi2023(LinCNS)} is dependent on the following proposition, which defines a time-dependent low-frequency cut-off function, $\hat{\Psi}$ in order to bound the kernel of the semigroup from below in the $L^\infty$-norm.

First, we denote 
\begin{align*}
    \x_{t} & \coloneqq (\x_1, t^{+1/4}\x_2, t^{+1/4}\x_3), \\
    \x_{t^{-1}} & \coloneqq (\x_1, t^{-1/4}\x_2, t^{-1/4}\x_3).
\end{align*}
We also take a nonnegative nonzero function $\hps \in C^\infty_0$ such that
\[
\supp \hps \subseteq \{ \x\in\R^3 \ | \ |\x|\in({1}/{2} , 1), \  |\x_1| \geq 1/2 \}, \quad
\hps(-\x) = \hps(\x), \text{ for all } \x \in \R^3.
\]

\begin{prop} {\rm (}\cite{ohaodha-iwabuchi2023(LinCNS)}{\rm )}
There exists a constant $C$ such that, for all $t$ sufficiently large, 
\begin{align*}
    \bignorm{ \F^{-1} \Bigf{ e^{t\lambda_\pm}
    \hps(t^{1/2}\x_{t}) 
    } }_\infty \geq Ct^{-2}.
\end{align*}
\end{prop}

\begin{rem}
We note that this bound from below on the low-frequency estimate is sufficient to prove that for all $t$ sufficiently large, 
\begin{align*}
    \bignorm{ \F^{-1} \Bigf{ \sum_{j\leq 2} \hp_j e^{t \lambda_\pm} } }_\infty \geq Ct^{-2}.
\end{align*}
Indeed, for all $t\geq 1$, we get by a simple application of Young's convolution inequality:
\begin{align*}
    \bignorm{ \F^{-1} \Bigf{ e^{t \lambda_\pm} 
    \hps(t^{1/2}\x_{t})
    } }_\infty
    & \leq
    \bignorm{ \F^{-1} \Bigf{ \hps(t^{1/2}\x_{t}) } }_1
    \bignorm{ \F^{-1} \Bigf{ \sum_{j\leq 2} \hp_j e^{t \lambda_\pm} } }_\infty
    \\
    & =
    \bignorm{ \F^{-1} \Bigf{ \hps(\x) } }_1
    \bignorm{ \F^{-1} \Bigf{ \sum_{j\leq 2} \hp_j e^{t \lambda_\pm} } }_\infty
    \\
    & \leq
    \norm{  \hps(\x) }_{W^{2,2}}
    \bignorm{ \F^{-1} \Bigf{ \sum_{j\leq 2} \hp_j e^{t \lambda_\pm} } }_\infty
    \\
    & \leq 
    C \bignorm{ \F^{-1} \Bigf{ \sum_{j\leq 2} \hp_j e^{t \lambda_\pm} } }_\infty.
\end{align*}
\end{rem}

We can extend this bound from below to the $\db^0_{ \infty , \infty }$-norm and thus obtain \eqref{p4 2nd main thm bnd from below} by proving the following proposition.

\begin{prop}\label{p4 Appendix bnd from below b0infin}
Let $\x_t$, $\hps$ be defined as above. Then for all $t$ sufficiently large, there exists a constant $C>0$ such that
\begin{align*}
    \bignorm{ 
    \F^{-1} \Bigf{ e^{t \lambda_\pm} } 
    }_{ \db^0_{\infty , \infty} }
    \geq C
    \bignorm{ \F^{-1} \Bigf{ e^{t \lambda_\pm} 
    \hps(t^{1/2}\x_{t})
    } }_\infty
    .
\end{align*}
\end{prop}

\begin{pf}
Note that, for all $t$ sufficiently large, 
\begin{align*}
    \supp{\hps} (t^{\frac{1}{2}} \x_t )
    & \subseteq
    \{ \x \in \R^3 \ | \ |\x_t| \in (  t^{-\frac{1}{2}} / 2 , \, t^{-\frac{1}{2}}) \}
    \\
    & \subseteq
    \{ \x \in \R^3 \ | \ |\x| \in (  t^{-\frac{1}{2}} / 2 , \, t^{-\frac{1}{2}} ) \}.
\end{align*}
Then, defining
\[
j_0 \coloneqq 1 - \floor{ \frac{1}{2} \log_2(t) },
\]
we get
\[
(  t^{-\frac{1}{2}} / 2 , \, t^{-\frac{1}{2}} )
\subseteq
(2^{j_0 - 2} , 2^{j_0 + 2}),
\]
and thus
\begin{align*}
    & 
    \bignorm{ \F^{-1} \Bigf{ e^{t \lambda_\pm} 
    \hps(t^{1/2}\x_{t})
    } }_\infty
    =
    \bignorm{ \F^{-1} \Bigf{ 
    \sum_{|j-j_0|\leq 2} \hat{\phi}_j
    e^{t \lambda_\pm} 
    \hps(t^{1/2}\x_{t})
    } }_\infty
    \\
    & \leq C
    \bignorm{ \F^{-1} \Bigf{ 
    \sum_{|j-j_0|\leq 2}\hat{\phi}_j
    e^{t \lambda_\pm} 
    } }_\infty
    \\
    & \leq C
    \bignorm{ \F^{-1} \Bigf{ 
    e^{t \lambda_\pm} 
    } }_{\db^0_{\infty , \infty} },
\end{align*}
for all sufficiently large $t>0$.
\end{pf}

\bigskip

\noindent {\bf Conflict of interest statement.} On behalf of all authors, the corresponding author states that there is no conflict of interest. 

\phantom{}

\noindent {\bf Data availability statement.} This manuscript has no associated data. 

\phantom{}

\bibliographystyle{unsrt} \bibliography{sample}

\end{document}